\documentclass[12pt]{amsart}%{article}%
\usepackage{amssymb}
\usepackage[mathscr]{eucal}
\usepackage{epsf}
\usepackage{epsfig}
\usepackage{color}
\DeclareGraphicsExtensions{.pstex,.eps,.epsi,.ps}

\makeatletter
\def\LaTeX{\leavevmode L\raise.42ex
    \hbox{\kern-.3em\size{\sf@size}{0pt}\selectfont A}\kern-.15em\TeX}
\makeatother

\newcommand{\BibTeX}{{\rm B\kern-.05em{\sc i\kern-.025emb}\kern-.08em\TeX}}

\newtheorem{thm}{Theorem}[section]
\newtheorem{prop}[thm]{Proposition}
\newtheorem{prop-def}[thm]{Proposition-Definition}
\newtheorem{lem}[thm]{Lemma}
\newtheorem{cor}[thm]{Corollary}
\newtheorem{conj}[thm]{Conjecture}

\newtheorem{defn}[thm]{Definition}
\newtheorem{exmp}[thm]{Example}

\newtheorem{rem}[thm]{Remark}

\newcommand{\proofbegin}{\noindent{\it Proof\,\,}}
\newcommand{\proofend}{\hfill$\Box$\bigskip}
\newenvironment{ackn}{\medskip \noindent \small
{\sl Acknowledgments.}}{\bigskip}

\makeatletter
\def\@currentlabel{2.1}\label{e:dispaa}
\def\@currentlabel{2.21}\label{e:dispau}
\def\@currentlabel{2.22}\label{e:dispav}
\def\@currentlabel{2.23}\label{e:dispaw}
\def\@currentlabel{2.24}\label{e:dispax}
\def\theequation{\thesection.\@arabic\c@equation}
\makeatother

\makeatletter
\def\alphenumi{%
  \def\theenumi{\alph{enumi}}%
  \def\p@enumi{\theenumi}%
  \def\labelenumi{(\@alph\c@enumi)}}
\makeatother

\newcommand{\e}{\epsilon}

\newcommand{\I}{\mathcal{I}}

\newcommand{\T}{\mathcal{T}}
\newcommand{\W}{\mathcal{W}}

\newcommand{\Real}{\mathbb{R}}

\newcommand{\Lg}{\mathfrak{g}}

\newcommand{\diff}{\mathcal{D}}
\newcommand{\sdiff}{\mathcal{D}_\mu}
\newcommand{\vect}{\mathfrak{X}}
\newcommand{\svect}{\mathfrak{X}_\mu}
\newcommand{\LD}{\mathcal{L}}
\newcommand{\td}{\tilde{d}}
\newcommand{\di}{\text{div}}
\newcommand{\ddt}{\frac{d}{dt}}

\newcommand{\pddt}{\frac{\partial}{\partial t}}

\newcommand{\ddtz}{\frac{d}{dt}\Big |_{t=0}}

\newcommand{\vol}{\text{vol}}

\newcommand{\Ent}{\text{Ent}}
\newcommand{\marginnote}[1]
{%\mbox{}\marginpar{\center{\hspace{0pt}\tiny{\bf#1}}}
}

\newcounter{bk}

\newcounter{pl}

\begin{document}
\bigskip
\bigskip

\title[A nonholonomic  Moser theorem and optimal transport]
{A nonholonomic Moser theorem \\and optimal mass transport}
\author{Boris Khesin}
%\address{Department of Mathematics, University of Toronto,
%ON M5S 2E4, Canada}
%\email{khesin@math.toronto.edu}
%\thanks{The author was partially supported by an NSERC grant.}
\author{Paul Lee}
\email{khesin@math.toronto.edu, plee@math.toronto.edu}
\address{Department of Mathematics, University of Toronto, ON M5S 2E4, Canada}

%\keywords{}
%  Math Subject Classifications
%\subjclass{}
\date{\today}
\maketitle

\begin{abstract}
We prove the following nonholonomic version of the classical Moser
theorem: given a bracket-generating distribution on a connected compact manifold
(possibly with boundary),
two volume forms of equal total volume can be
isotoped by the flow of a vector field tangent to this distribution.
We describe formal solutions of the corresponding nonholonomic mass transport
problem and present the Hamiltonian framework for both the Otto calculus and its
nonholonomic counterpart as infinite-dimensional Hamiltonian reductions
on diffeomorphism groups.

Finally, we define a nonholonomic analog of the Wasserstein (or, Kantorovich) metric
on the space of densities  and prove that the subriemannian
heat equation defines a gradient flow on the nonholonomic Wasserstein space
with the potential given by the Boltzmann relative entropy functional.
\end{abstract}

\tableofcontents
%%% ----------------------------------------------------------------------

\section{Introduction}

The classical Moser theorem establishes that the total volume is the only
invariant for a volume form on a compact connected manifold with respect to
the diffeomorphism action. In this paper we prove a nonholonomic counterpart
of this result and present its applications in the problems of nonholonomic
optimal mass transport.

The equivalence for the diffeomorphism action is often formulated in terms of
``stability" of the corresponding object: the existence of a
diffeomorphism relating the initial object with a deformed one means
that the initial object is stable, as it differs from the deformed
one merely by a coordinate change. Gray showed in \cite{Gray1959}
that contact structures on a compact manifold are stable.
Moser \cite{Moser1965} established stability for volume forms and
symplectic structures. A leafwise counterpart of Moser's argument for foliations
was presented by Ghys in \cite {Ghys1984},
while stability of symplectic-contact pairs in
transversal foliations was proved in
\cite{Bande-Ghiggini-Kotschick2004}. In this paper we
establish stability of volume forms in the presence of any
bracket-generating distributions on connected compact manifolds:  two volume forms
of equal total volume on such a manifold can be
isotoped by the flow of a vector field tangent to the distribution.
We call this statement a nonholonomic  Moser theorem.

%Polterovich 2000: For a contact structure, contact diffeomorphisms have an
%invariant w.r.t. a volume form, i.e. pairs (diffeomorphism, volume form)
%are not stable. (Note: contact diffeomorphisms are NOT time-one-maps of
%contact vector fields.)
\smallskip

Recall that a distribution $\tau$ on the manifold $M$ is called {\it bracket-generating},
or {\it completely nonholonomic}, if local vector fields
tangent to $\tau$ and their iterated Lie brackets span the entire
tangent bundle of the manifold $M$.
Nonholonomic distributions arise in various problems related to
rolling or skating, wherever the ``no-slip" condition is present. For
instance, a ball rolling over a table defines a trajectory in a
configuration space tangent to a nonholonomic distribution of
admissible velocities. Note that such a ball can be rolled to any
point of the table and stopped at any \textit{a priori} prescribed position.
The latter is a manifestation of the Chow-Rashevsky theorem (see e.g. \cite{Mo}):
For a bracket-generating  distribution $\tau$ on a connected
manifold $M$ any two points in $M$ can be connected by a horizontal
path (i.e. a path everywhere tangent to the distribution $\tau$).\footnote{
The motivation for considering volume forms (or, densities) in a space with
distribution can be related to problems with many tiny rolling balls:
It is more convenient to consider the density of such balls, rather
than look at them individually.}

Note that for an integrable distribution there is a foliation
to which it is tangent and a horizontal path always stays on
the same leaf of this foliation. Furthermore, for an integrable distribution,
the existence of an isotopy between volume forms requires an infinite
number of conditions. On the contrary, the nonholonomic Moser theorem
shows that a non-integrable bracket-generating
distribution imposes only one condition on total volume of the
forms for the existence of the isotopy between them.

Closely related to the nonholonomic Moser theorem
is the existence of a nonholonomic Hodge decomposition, and
the corresponding properties of the
subriemannian Laplace operator, see Section \ref{nonHodge}.
We also formulate the corresponding nonholonomic mass transport
problem and describe its formal solutions as projections of horizontal
geodesics on the diffeomorphism group for the $L^2$-Carnot-Caratheodory
metric.

In order to give this description, we first present
the Hamiltonian framework for what is now called the Otto calculus -
the Riemannian submersion picture for the problems of optimal mass transport.
It turns out that the submersion properties can be naturally understood
as an infinite-dimensional Hamiltonian reduction
on diffeomorphism groups, and this admits a generalization to the
nonholonomic setting. We define a nonholonomic analog of the Wasserstein metric
on the space of densities. Finally, we extend Otto's result on the heat equation
and  prove that the subriemannian
heat equation defines a gradient flow on the nonholonomic Wasserstein space
with potential given by  the Boltzmann relative entropy functional.

\bigskip

%%%%%%%%%%%%%%%%%%%%%%%%%%%%%%%%%
%%%%%%%%%%%%%%%%%%%%%%%%%%%%%%%%%

\section{Around Moser's theorem}\label{around}

\medskip

\subsection{Classical and nonholonomic Moser theorems}
The main goal of this section is to prove the following
nonholonomic version of the classical Moser theorem. Consider  a
distribution $\tau$ on a compact manifold $M$
(without boundary unless otherwise stated).

\begin{thm}\label{thm:nonhol_vol}
Let $\tau$ be a bracket-generating  distribution,
and $\mu_0, \,\mu_1$ be two volume forms on $M$
with the same total volume: $\int_M  \mu_0=\int_M  \mu_1$.
 Then there exists
a diffeomorphism $\phi$ of $M$ which is the time-one-map of the flow $\phi_t$
of a non-autonomous vector field $V_t$ tangent to the distribution
$\tau$ everywhere on $M$  for every $t\in [0,1]$, such that
$\phi^*\mu_1=\mu_0$.
\end{thm}

Note that the existence of the ``nonholonomic isotopy" $\phi_t$ is
guaranteed by the only condition on equality of total volumes
for  $\mu_0$ and $\mu_1$, just like in the classical case:

\begin{thm}\label{thm:clas_mos}\cite{Moser1965}
Let $M$ be a manifold without boundary, and $\mu_0, \,\mu_1$ are two
volume forms on $M$ with the same total volume: $\int_M \mu_0=\int_M
\mu_1$. Then there exists a diffeomorphism $\phi$ of $M$, isotopic
to the identity, such that $\phi^*\mu_1=\mu_0$.
\end{thm}
\medskip

\begin{rem} \label{re:clas_mos}
{\rm
The classical Moser theorem has numerous variations and generalizations,
some of which we would like to mention.

a) Similarly one can show that not only the identity, but
any diffeomorphism of $M$ is isotopic
to a diffeomorphism which pulls back $\mu_1$ to $\mu_0$.

b) The Moser theorem also holds for a manifold $M$ with boundary.
In this case a diffeomorphism $\phi$ is a time-one-map for a (non-autonomous)
vector field $V$ on $M$, tangent to the boundary $\partial M$.

c)  Moser also proved in \cite{Moser1965}
a similar statement for a pair of symplectic forms on a
manifold $M$: if two symplectic structures can be deformed to each other
among symplectic structures in the same cohomology class on $M$,
these deformation can be carried out by a flow of diffeomorphisms of $M$.

Below we describe to which degree
these variations extend to the nonholonomic case.
}
\end{rem}

\medskip
%%%%%%%%%%%%%%%%%%%%%%%%%%%%%%%%%%%

\subsection{The Moser theorem for a fibration}

Apparently, the most straightforward generalization of the classical
Moser theorem is its version ``with parameters." In this case, volume
forms on $M$ smoothly depend on parameters and
have the same total volume at each value of this
parameter: $\int_M  \mu_0(s)=\int_M  \mu_1(s)$ for all $s$.
The theorem guarantees that the corresponding
diffeomorphism exists and depends smoothly on this parameter $s$.

The following theorem can be regarded as a modification
of the parameter version:

\begin{thm}\label{thm:fibration}
Let $\pi:N\to B$ be a fibration of an $n$-dimensional manifold $N$
over a $k$-dimensional base manifold $B$. Suppose that $\mu_0,
\,\mu_1$ are two smooth volume forms on $N$. Assume that the
pushforwards of these $n$-forms to $B$ coincide, i.e. they give one
and the same $k$-form on $B$: $\pi_*\mu_0=\pi_*\mu_1$. Then, there
exists a diffeomorphism $\phi$ of $N$ which is the time-one-map of a
(non-autonomous) vector field $V$ tangent everywhere to the fibers
of this fibration and such that $\phi^*\mu_1=\mu_0$.
\end{thm}

\begin{rem}
{\rm
Note that in this version
the volume forms are given on the {\it ambient} manifold $N$, while
in the parametric version of the Moser theorem we are given
{\it fiberwise} volume forms. There is also
a similar version of this theorem for a foliation, cf. e.g.
\cite{Ghys1984}. In either case, for the corresponding diffeomorphism
to exist the volume forms have to satisfy infinitely many conditions
(the equality of the total volumes as functions in the parameter $s$ or as
the push-forwards $\pi_*\mu_0$ and $\pi_*\mu_1$).
The case of a fibration (or a foliation) corresponds to an
integrable distribution $\tau$, and presents the ``opposite case"
to a bracket-generating distribution.
Unlike the case of an integrable distribution, the existence of
the corresponding isotopy between volume forms in the bracket-generating case
imposes only one condition, the equality of the total volumes of the two forms
(regardless, e.g., of the distribution growth vector  at different points
of the manifold).
}
\end{rem}

\medskip

%%%%%%%%%%%%%%%%%%%%%%%%%%%%%%%%%%%

\subsection{Proofs}

First, we recall a proof of the classical Moser theorem. To show how
the proof changes in the nonholonomic case, we split it into several steps.

\begin{proof}
1) Connect the volume forms $\mu_0$ and $\mu_1$ by a ``segment''
$\mu_t=\mu_0+t(\mu_1-\mu_0),$ $t\in[0,1]$. We will be looking for
a diffeomorphism $g_t$ sending $\mu_t$ to $\mu_0$:
$g_t^*\mu_t=\mu_0$.
%Let $v_t$ be the velocity of the 1-parameter family of diffeomorphisms $g_t$:
By taking the $t$-derivative of this  equation, we get the following
``homological equation" on the velocity $V_t$ of the flow $g_t$:
$g_t^*(\LD_{V_t}\mu_t+\partial_t\mu_t)=0$, where $\partial_t
g_t(x)=V_t(g_t(x))$. This is equivalent to
$$
\LD_{V_t}\mu_t=\mu_0-\mu_1\,,
$$
since $\partial_t\mu_t=-(\mu_0-\mu_1)$.

By rewriting $\mu_0-\mu_1=\rho_t\mu_t$ for an appropriate function $\rho_t$,
we reformulate the equation
$\LD_{V_t}\mu_t=\rho_t\mu_t$ as the problem $\di_{\mu_t}V_t=\rho_t$
of looking for a vector field $V_t$ with a prescribed divergence $\rho_t$.
Note that the total integral of the
function $\rho_t$ (relative to the volume $\mu_t$) over $M$ vanishes, which manifests
the equality of total volumes for $\mu_t$.

2) We omit the index $t$ for now and consider a Riemannian metric on
$M$ whose volume form is $\mu$. We are looking for a required
field $V$ with prescribed divergence among gradient vector fields
$V=\nabla u$, which ``transport the mass" in the fastest way.
This leads us to the elliptic equation $\di_\mu(\nabla u)=\rho$, i.e.
 $\Delta u=\rho$, where the Laplacian $\Delta$ is
defined by $\Delta u:=\di_\mu\nabla u$ and depends on the
Riemannian metric on $M$.

3) The key part of the proof is the following

\begin{lem}\label{elliptic}
The Poisson equation $\Delta u=\rho$ on a compact Riemannian
manifold $M$ is solvable for
any function $\rho$ with zero mean: $\int_M\rho\,\mu=0$
(with respect to the Riemannian volume form $\mu$).
\end{lem}

\proofbegin of Lemma. Describe the space $\text{Coker}\,
\Delta:=(\text{Im}\, \Delta)^{\perp_{L^2}}$, i.e. find
the space of all functions $h$ which are
$L^2$-orthogonal to the image $\text{Im}\,\Delta$. By applying
integration by parts twice, one has:
$$
0=\langle h,\Delta u\rangle _{L^2}=-\langle \nabla h,\nabla u\rangle _{L^2} =
\langle \Delta h, u\rangle _{L^2}
$$
for all smooth functions $u$ on $M$. Then such functions $h$ must be
harmonic, and hence they are constant functions on $M$:
$(\text{Im}\, \Delta)^{\perp_{L^2}}=\{const\}$.
Since the image $\text{Im}\, \Delta$ is closed, it is the $L^2$-orthogonal
complement of the space of constant functions
$\text{Im}\, \Delta=\{const\})^{\perp_{L^2}}$.
The condition of orthogonality to
constants is exactly the condition of zero mean for $\rho$:
$\langle const,\rho\rangle _{L^2}=\int_M \rho\,\mu=0$. Thus the equation  $\Delta u=\rho$
has a weak solution for $\rho$ with zero mean, and
the ellipticity of $\Delta$ implies that the solution is smooth for a smooth
function $\rho$.
\proofend

4) Now, take $V_t:=\nabla u_t$ and let $g^t_{V}$ be the
corresponding flow on $M$. Since $M$ is compact and $V_t$ is smooth,
the flow exists for all time $t$. The diffeomorphism
$\phi:=g^1_{V}$, the time-one-map of the flow $g^t_{V}$, gives the
required map which pulls back the volume form $\mu_1$ to $\mu_0$:
$\phi^*\mu_1=\mu_0$.
\end{proof}
\bigskip

\proofbegin of Theorem \ref{thm:fibration}, the Moser theorem for a fibration:

We start by defining the new volume form on the fibres $F$ using the
pushforward $k$-form $\nu_0:=\pi_*\mu_0$ on the base $B$ and the
volume $n$-form $\mu_0$ on $N$. Namely, consider the pull-back
$k$-form $\pi^*\nu_0$ to $N$. Then there is a unique $(n-k)$-form
$\mu_0^F$ on fibers such that $\mu_0^F\wedge \pi^*\nu_0=\mu_0$.
Similarly we find $\mu_1^F$. Due to the equality of the pushforwards
$\pi_*\mu_0$ and $\pi_*\mu_1$, the total volumes of $\mu_0^F$ and
$\mu_1^F$ are fiberwise equal. Hence by the Moser theorem applied to
the fibres, there is a smooth vector field tangent to the fibers,
smoothly depending on a base point, and  whose flow sends
one of the $(n-k)$-forms, $\mu_1^F$, to the other,
$\mu_0^F$. This field is defined globally on $N$, and hence its
time-one-map pulls back $\mu_1$ to $\mu_0$.
\proofend
\bigskip

Now we turn to a  nonholonomic distribution on a manifold.

\medskip

\proofbegin of Theorem \ref{thm:nonhol_vol},
the nonholonomic version of the Moser theorem.

1) As before, we connect the forms by a segment $\mu_t, t\in[0,1]$,
and we come to the same homological equation. The latter reduces to
$\di_\mu V=\rho$ with $\int\rho\,\mu=0$, but the equation now
is for a vector field $V$ tangent to the distribution $\tau$.

2) Consider some Riemannian metric on $M$. Now we will be looking
for the required field $V$ in the form $V:=P^\tau\nabla u$,
where $P^\tau$ is a pointwise orthogonal projection of tangent vectors
to the planes of our distribution $\tau$.

We obtain the equation $\di_\mu(P^\tau\nabla u)=\rho$.
Rewrite this equation by introducing the sub-Laplacian
$\Delta^\tau u:=\di_\mu(P^\tau\nabla u)$ associated to the
distribution $\tau$ and the {\it Riemannian} metric on $M$. The equation
on the potential $u$  becomes $\Delta^\tau u=\rho$.

3) An analog of Lemma \ref{elliptic} is now as follows.

\begin{prop}\label{prop:key-prop}
a) The sub-Laplacian operator $\Delta^\tau u:={\rm div }_\mu
(P^\tau\nabla u)$ is a self-adjoint hypoelliptic operator. Its
image is closed in $L^2$.

b) The equation $\Delta^\tau u=\rho$ on a compact Riemannian
manifold $M$ is solvable for any function $\rho$ with zero mean:
$\int_M\rho\,\mu=0$.
\end{prop}

\proofbegin of Proposition. a) The principal symbol
$\delta^\tau$ of the operator $\Delta^\tau$ is the sum of
squares of vector fields forming a basis for the
distribution $\tau$: $\delta^\tau=\sum X_i^2$, where $X_i$ form a
horizontal orthonormal frame for $\tau$.
This is exactly the H\"ormander condition of
hypoellipticity \cite{Hor1} for the operator $\Delta^\tau$.
The self-adjointness follows from the properties of projection and
integration by parts. The closedness of the image in $L^2$ follows
from the results of \cite{RoSt, RoTar}.

b) We need to find the condition of weak solvability in $L^2$ for
the equation $\Delta^\tau u=\rho$. Again, we are looking for all those 
functions $h$ which are $L^2$-orthogonal to the image of $\Delta^\tau$ (or, which
is the same, in the kernel of this operator):
$$
0=\langle h,\Delta^\tau u\rangle _{L^2}=\langle h,{\rm div }_\mu
(P^\tau\nabla u)\rangle _{L^2}
$$
for all smooth functions $u$ on $M$. In particular, this should hold for $u=h$.
Integrating by parts we come to
$$
0=\langle h,{\rm div }_\mu (P^\tau\nabla h)\rangle _{L^2}
=-\langle \nabla h, P^\tau\nabla h\rangle _{L^2}
=-\langle P^\tau\nabla h, P^\tau\nabla h\rangle _{L^2}\,,
$$
where in the last equality we used the projection property $(P^\tau)^2=P^\tau
=(P^\tau)^*$.
Then $P^\tau\nabla h=0$ on $M$, and hence the equation $\Delta^\tau u=\rho$ is
solvable for any function $\rho\perp_{L^2}\{h~|~P^\tau\nabla h=0\}$.
We claim that all such functions $h$
are constant on $M$. Indeed, the condition $P^\tau\nabla h=0$
means that $\LD_{X}h=0$
for any horizontal field $X$, i.e. a field tangent to the distribution $\tau$.
But then $h$ must be constant
along any horizontal path, and due to the Chow-Rashevsky theorem it
must be constant everywhere on $M$. Thus the functions $\rho$ must
be $L^2$-orthogonal to all constants, and hence they have zero
mean. This implies that the equation $\di_\mu(P^\tau\nabla u)=\rho$ is solvable
for any $L^2$ function $\rho$ with zero mean. For a smooth $\rho$ the solution
is also smooth due to hypoellipticity of the operator.
\proofend

4) Now consider the horizontal field $V_t:=P^\tau\nabla u_t$.
As before, the time-one-map of its flow exists for the smooth field $V_t$ on the
compact manifold $M$, and it  gives the required diffeomorphism $\phi$.
\proofend

\bigskip
%%%%%%%%%%%%%%%%%%%%%%%%%%%%%%%%%%%

\subsection{The nonholonomic Hodge decomposition and sub-Laplacian}\label{nonHodge}
\begin{figure}[ht!]
\input epsf
\centerline{\epsfysize=0.33\vsize\epsffile{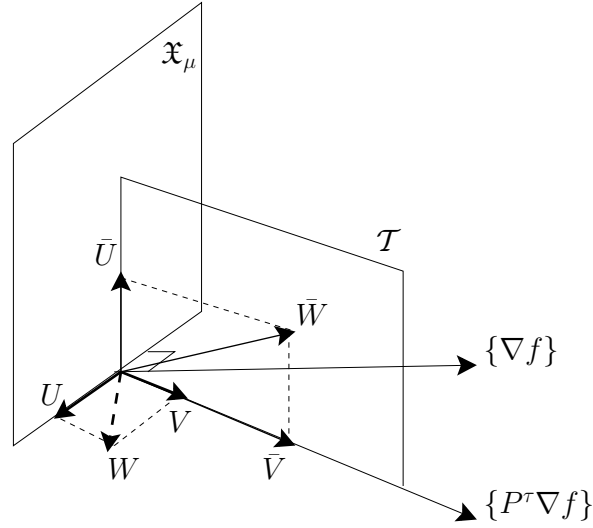}}
 \caption{A nonholonomic Hodge decomposition.}\label{fig:hodge}
\end{figure}

According to the classical Helmholtz-Hodge decomposition, any vector
field $W$ on a Riemannian manifold $M$ can be uniquely decomposed
into the sum $W=\tilde V+\tilde U$, where $\tilde V=\nabla f$ and $\di_\mu \tilde U=0$.
Proposition \ref{prop:key-prop} suggests the following nonholonomic Hodge
decomposition of vector fields on a manifold with a
bracket-generating distribution:

\begin{prop}\label{prop:non-hol-hodge}
1) For a bracket-generating distribution $\tau$ on a Riemannian
manifold $M$, any vector field $W$ on $M$  can be uniquely
decomposed into the sum $W=V+U$, where the field $V=P^\tau\nabla f$
and it is tangent to the distribution $\tau$, while the field
$U$ is divergence-free: ${\rm div }_\mu\,U=0$. Here $P^\tau$ is the
pointwise orthogonal projection to $\tau$.

2) Moreover, if the vector field $\bar W$ is tangent to the
distribution $\tau$ on $M$, then $ \bar W=\bar V+\bar U$, where
$\bar V=P^\tau\nabla f\,||\, \tau$ as before, while the field $\bar U$
is divergence-free, tangent to $\tau$, and $L^2$-orthogonal to $\bar V$,
see Figure \ref{fig:hodge}.
\end{prop}
\bigskip

%\begin{figure}[ht!]
%\begin{center}
%\centerline{\scalebox{0.5}{\input{mostranfig1.pstex_t}}}
% \caption{A nonholonomic Hodge decomposition.}\label{fig:hodge}
%\end{center}
%\end{figure}

%\bigskip

\begin{proof}
Let $\rho:=\di_\mu\,W$ be the divergence of $W$ with respect to the
Riemannian volume $\mu$. First, note that $\int_M\rho\,\mu=0$.
Indeed, $\int_M (\di_\mu\,W)\,\mu=\int_M \LD_W\mu=0$, since the volume
of $\mu$ is defined in a coordinate-free way, and does not change
along the flow of the field $W$.

Now, apply Proposition \ref{prop:key-prop} to find a solution of the
equation $\di\,(P^\tau\nabla f)=\rho$. The field
$V:=P^\tau\nabla f$ is defined uniquely. Then the field $U:=W-V$
is divergence-free, which proves $1)$.

For a field $\bar W\,||\,\tau$, we define $\bar V:=P^\tau\nabla f$
in the same way. Note that $\bar V||\,\tau$ as well. Then  $\bar U:=\bar
W-\bar V$ is both tangent to $\tau $ and divergence-free.
Furthermore,
$$
\langle \bar U, V\rangle _{L^2}=\langle \bar U, P^\tau\nabla f\rangle _{L^2}=\langle P^\tau\bar U,
\nabla f\rangle _{L^2} =\langle \bar U, \nabla f\rangle _{L^2}=\langle  {\rm div }_\mu\,\bar U,
f\rangle _{L^2}=0,
$$
where we used the properties of $\bar U$ established above: $\bar
U\,||\,\tau$ and $\di_\mu\,\bar U=0$.

\end{proof}

Above we defined a sub-Laplacian $\Delta^\tau u:=\di_\mu\,(P^\tau\nabla u)$
for a function $u$
on a Riemannian manifold $M$ with a distribution $\tau$.

\begin{prop}(cf. \cite{Mo})
The sub-Laplacian $\Delta^\tau $
depends only on a subriemannian metric on the distribution
$\tau$ and a volume form in the ambient manifold $M$.
\end{prop}

\begin{proof}
Note that the operator $P^\tau\nabla$ on a function $u$ is
the horizontal gradient $\nabla^\tau$ of $u$, i.e.
the vector of the fastest growth of $u$ among the directions in
$\tau$. If one chooses a local orthonormal
frame $X_1,...,X_k$ in $\tau$, then $P^\tau\nabla u=
\sum_{i=1}^k(\LD_{X_i}\,u)X_i$. Thus
the definition of the horizontal gradient relies on the
subriemannian metric only.

The sub-Laplacian $\Delta^\tau\psi=\di_\mu\,(P^\tau\nabla\psi)$
needs also the volume form
$\mu$ in the ambient manifold to take the divergence with respect to
this form.
\end{proof}

The corresponding {\it nonholonomic
heat equation} $\partial_t u=\Delta^\tau u$ is also defined
by the subriemannian metric and a volume form.

%In the case  of a contact distribution $\tau$, appearing in $CR$ structures, the corresponding asymptotics have been studied in \cite{Beals and K}. Apparently, the asymptotic in the general case have to depend, in particular, on the growth vector of $\tau$.
\medskip

%%%%%%%%%%%%%%%%%%%%%%%%%%%%%%%%%%%%%%%%%%%%%%%%%%%%%%%%

\bigskip

\subsection{The case with boundary}
For a manifold $M$ with non-empty boundary $\partial M$ and two volume forms
$\mu_0,\mu_1$ of equal total volume, the
classical Moser theorem establishes the existence of diffeomorphism $\phi$
which is the time-one-map for the flow of a field $V_t$ tangent to $\partial M$
and such that $\phi^*\mu_1=\mu_0$.

The existence of the required gradient field $V_t=\nabla u$ is guaranteed
by the following

\begin{lem}\label{elliptic_bound}
Let $\mu$ be a volume form on a Riemannian manifold $M$ with boundary $\partial M$.
The Poisson equation $\Delta u=\rho$ with Neumann boundary condition
$\frac{\partial}{\partial n} u=0$ on the boundary $\partial M$
is solvable for any function $\rho$ with zero mean: $\int_M\rho\,\mu=0$.
\end{lem}

Here $\frac{\partial}{\partial n}$ is the differentiation in the direction of outer normal
$n$ on the boundary.

\medskip
\proofbegin of Lemma.
Proceed in the same way as in Lemma \ref{elliptic} to find all functions
$L^2$-orthogonal to the image $\text{Im}\,\Delta$. The first
integration by parts gives:

$$
0=\int_M h(\Delta u)\mu
=-\int_M\langle \nabla h,\nabla u\rangle \mu
           +\int_{\partial M}h\,(\frac{\partial}{\partial n}u)\,\mu
=-\int_M\langle \nabla h,\nabla u\rangle \mu\,,
$$
where in the last equality we used the Neumann boundary conditions.
The second integration by parts gives:
$$
0=\int_M\langle \Delta h, u\rangle \mu
-\int_{\partial M}(\frac{\partial}{\partial n}h)\,u\,\mu
$$
This equation holds for all smooth functions $u$ on $M$,
so any such function $h$ must be harmonic in $M$ and satisfy the Neumann
boundary condition $\frac{\partial}{\partial n}h=0$. Hence, these
are constant functions on $M$: $(\text{Im}\, \Delta)^{\perp_{L^2}}=\{const\}$.
This gives the same description as in the no-boundary case:
the image $(\text{Im}\, \Delta)$ with the Neumann  condition
consists of functions $\rho$ with  zero mean.
\proofend

Geometrically, the Neumann boundary condition means that there is no flux of
density through the boundary $\partial M$: $0=\frac{\partial u}{\partial n}
=n\cdot\nabla u=n\cdot V$ on $\partial M$.

\medskip

For distributions on manifolds with boundary, the solution of the Neumann
problem becomes a much more subtle issue, as the behavior of the distribution near
the boundary affects the flux of horizontal fields across the boundary, and hence
the solvability in this problem. However, there is a class of domains in
length spaces for which the solvability of the Neumann problem was established.

Let $LS$ be a length space with the distance function $d(x,y)$,  defined
as infimum of lengths of continuous curves joining $x,y\in LS$. Consider
domains in this space with the property that  sufficiently close points in
those domains can be joined by a
not very long path which does not get too close to the domain boundary. The formal
definition is as follows.

\begin{defn}
{\rm
An open set $\Omega\subseteq LS$ is called an {\it $(\e,\delta)$-domain} if there
exist $\delta>0$ and $0<\e\leq 1$ such that for any pair of points $p,q\in
\Omega$ with $d(p,q)\leq\delta$ there is a continuous rectifiable curve
$\gamma:[0,T]\to\Omega$ starting at $p$ and ending at $q$ such that
 the length $l(\gamma)$ of the curve $\gamma$ satisfies
$$
l(\gamma)\leq\frac{1}{\e}d(p,q)
$$
and
$$
\text{min}\{d(p,z),d(q,z)\}\leq\frac{1}{\e} d(z,\partial\Omega)
$$
for all points $z$ on the curve $\gamma$.
%The {\it diameter $\text{rad}(\Omega)$} of a domain $\Omega$ is
%$$
%diam(\Omega)=\sup\{r>0~|~\text{ for every } 0\leq s<r \text{ and any } p\in \Omega,
%\text{ there exists }q\in \Omega \text{ with }d(p,q)=s\}\,.
%$$
}
\end{defn}

A large source of $(\e,\delta)$-domains is given by some classes of
open sets in Carnot groups,
where the Carnot group itself is regarded as a length space with the
Carnot-Caratheodory distance, defined via the lengths of admissible (i.e. horizontal)
paths, see e.g. \cite{Nhieu01}. There is a natural notion of diameter (or, radius)
for domains in length spaces.

\begin{thm}\label{thm:nonhol_vol_bound}
Let $\tau$ be a bracket-generating  distribution
on a subriemannian manifold $M$ with smooth boundary $\partial M$,
and $\mu_0, \,\mu_1$ be two volume forms on $M$
with the same total volume: $\int_M  \mu_0=\int_M  \mu_1$.
Suppose that the interior of $M$ is an $(\e,\delta)$-domain of positive diameter.

Then there exists
a diffeomorphism $\phi$ of $M$ which is the time-one-map of the flow $\phi_t$
of a non-autonomous vector field $V_t$ tangent to the distribution
$\tau$ everywhere on $M$ and to the boundary $\partial M$
for every $t\in [0,1]$, such that $\phi^*\mu_1=\mu_0$.
\end{thm}

The proof immediately follows from the result on solvability of the corresponding
Neumann problem $\Delta^\tau u=\rho$ with $\,n\cdot ( P^\tau\nabla u)|_{\partial M}=0$
(or, which is the same, $\frac{\partial u}{\partial (P^\tau n)}|_{\partial M}=0$)
for such domains, established in
\cite{Nhieu, Nhieu01} (cf. Theorem 1.5 in \cite{Nhieu01}).

\bigskip

%%%%%%%%%%%%%%%%%%%%%%%%%%%%%%%%%%%

%%%%%%%%%%%%%%%%%%%%%%%%%%%%%%%%%%%

\section{Distributions on diffeomorphism groups}\label{sect: distributions}

\subsection{A fibration on the group of  diffeomorphisms}\label{subsec:fibra}

Let $\mathcal D$ be the group of all (orientation-preserving)
diffeomorphisms of a manifold $M$. Its Lie
algebra $\vect$ consists of all smooth
vector fields on  $M$. The tangent space to the
diffeomorphism group at any point $\phi\in\mathcal D$ is given by
the right translation of the Lie algebra $\vect$ from the identity
$id\in\diff$ to $\phi$:
\[
T_\phi\diff=\{X\circ\phi~|~X\in\vect\}\,.
\]

Fix a volume form $\mu$  of total volume 1
on the manifold $M$. Denote by $\diff_\mu$ the subgroup of volume-preserving
diffeomorphisms, i.e. the diffeomorphisms preserving the volume form $\mu$.
The corresponding Lie algebra $\vect_\mu$
is the space  of all vector fields on the manifold $M$ which are
divergence-free with respect to the volume form $\mu$.

\smallskip

Let $\W$ be the set of all smooth normalized volume forms in $M$,
which  is called the {\it (smooth) Wasserstein space.} Consider the
projection map $\pi^\diff:\diff\to\W$ defined by the push forward of
the fixed volume form $\mu$ by the diffeomorphism $\phi$, i.e.
$\pi^\diff(\phi)=\phi_*\mu$.
The projection $\pi^\diff:\diff\to\W$ defines a natural structure of
a principal bundle on $\diff$ whose structure group is the subgroup
$\diff_\mu$ of volume-preserving diffeomorphisms of $M$ and fibers
$F$ are right cosets for this subgroup in $\diff$. Two
diffeomorphisms $\phi$ and $\tilde \phi$ lie in the same fiber if they
differ by a composition (on the right) with a volume-preserving
diffeomorphism: $\tilde \phi=\phi\circ s, \,\,s\in\diff_\mu$.

On the  group $\diff$ we define two vector bundles $V\!er$ and $H\!or$ whose spaces
at a diffeomorphism $\phi\in \diff$ consist of
right translated divergence-free fields
$$
V\!er_\phi=\{X\circ\phi~|~\di_{\phi_*\mu} X=0\}
$$
and gradient fields
$$
H\!or_\phi:=\{\nabla f\circ\phi~|~f\in C^\infty(M)\}\,,
$$
respectively. Note that the bundle $V\!er$ is defined by the fixed
volume form $\mu$, while $Hor$ requires a Riemannian metric.\footnote{The metric
on $M$ does not need to have the volume form $\mu$.
In the general case, $\svect$ consists of vector fields divergence-free with
respect to $\mu$, while the gradients are considered for the chosen metric on $M$.}

\begin{prop}
The bundle $V\!er$ of translated divergence-free fields is the bundle of vertical
spaces $T_\phi F$ for the fibration $\pi^\diff:\diff\to\W$.
The bundle $H\!or$ over $\diff$ defines a
horizontal distribution for this fibration $\pi^\diff$.
\end{prop}

\begin{proof}
Let $\phi_t$ be a curve in a fibre of $\pi^\diff:\diff\to\W$
emanating from the point $\phi_0=\phi$. Then $\phi_t=\phi_0\circ s_t$,
where $s_0=id$ and $s_t$ are volume-preserving diffeomorphisms for each $t$.
Let $X_t$ be a family of divergence-free vector fields,
such that $\partial_t s_t=X_t\circ s_t$. Then the vector tangent to the curve
$\phi_t=\phi_0\circ s_t$ is given by
$\ddtz(\phi_0\circ s_t)=(\phi_{0*}X_0)\circ\phi_0$. Since $X_0$ is
divergence-free with respect to $\mu$, $\phi_{0*}X_0$ is divergence-free
with respect to $\phi_*\mu$. Hence, any vector tangent to the
diffeomorphism group at $\phi$ is given by $X\circ\phi$, where $X$
is a divergence-free field with respect to the form $\phi_*\mu$.

By the Hodge
decomposition of vector fields, we have the direct sum
$T\diff=H\!or\oplus V\!er$.
\end{proof}

\begin{rem}\label{pbdle}
{\rm
The classical Moser theorem
\ref{thm:clas_mos} can be thought of as the existence of path-lifting property
for the principal bundle $\pi^\diff:\diff\to\W$: any deformation
of volume forms can be traced by the corresponding flow,
i.e. a path on the diffeomorphism group, projected to the deformation of forms.
Its proof shows that this path lifting property holds and has the uniqueness
property in the presence of the horizontal distribution defined above by
using the Hodge decomposition. Namely, given
any path $\mu_t$ starting at $\mu_0$
in the smooth Wasserstein space $\W$ and a point
$\phi_0$ in the fibre $(\pi^\diff)^{-1}\mu_0$, there exists a unique
path $\phi_t$ in the diffeomorphism group which is tangent to the
horizontal bundle $Hor$, starts at $\phi_0$, and projects to $\mu_t$,
see Figure \ref{fig:lifting}.
\begin{figure}[ht!]
\input epsf
\centerline{\epsfysize=0.44\vsize\epsffile{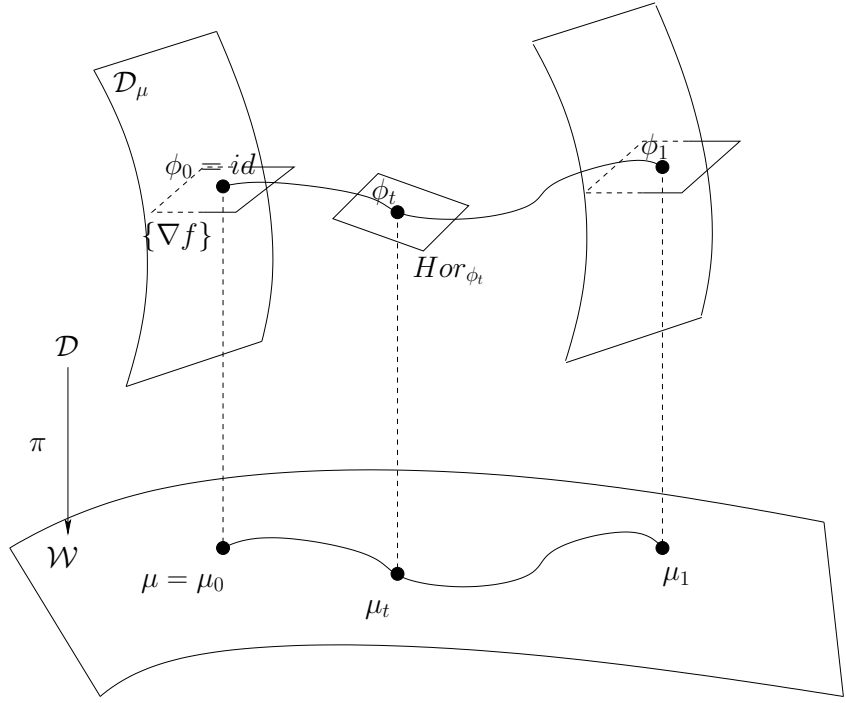}}
 \caption{The Moser theorem in both the classical and nonholonomic settings
is a path-lifting property in the diffeomorphism group.}\label{fig:lifting}
\end{figure}
}
\end{rem}
%\medskip

%\begin{figure}[ht!]
%\centerline{\scalebox{0.5}{\input{mostranfig2.pstex_t}}}
% \caption{The Moser theorem in both the classical and nonholonomic settings
%is a path-lifting property in
%the diffeomorphism group.}\label{fig:lifting}
%\end{figure}

\bigskip

%%%%%%%%%%%%%%%%%%%%%%%

\subsection{A nonholonomic distribution on the diffeomorphism group}

Let $\tau$ be a bracket-generating distribution on the manifold $M$.
Consider the right-invariant distribution $\mathcal T$ on the
diffeomorphism group $\diff$ defined at the identity $id\in \diff$
of the group by the subspace in $\vect$ of all those vector fields
which are tangent to the distribution $\tau$ everywhere on $M$:
$$
\mathcal T_\phi=\{V\circ\phi~|~V(x)\in\tau_x \,\text{ for all }\,x\in M\}.
$$
\medskip

\begin{prop}  The infinite-dimensional  distribution $\mathcal T$
is a non-integrable distribution in  $\diff$.
Horizontal paths in this distribution
are flows of non-autonomous vector fields tangent to the distribution
$\tau$ on manifold $M$.
\end{prop}

\begin{proof} To see  that the distribution $\mathcal T$
is non-integrable we consider two horizontal vector fields $V$ and $W$ on $M$
and the corresponding right-invariant vector
fields $\widetilde V$ and $\widetilde W$
on $\diff$. Then their bracket at the identity of the group
is (minus) their commutator as vector fields $V$ and $W$ in $M$.
This commutator does not belong to the plane ${\mathcal T}_{id}$
since the distribution $\tau$ is non-integrable,
and at least somewhere on $M$ the commutator of horizontal fields $V$ and $W$
is not horizontal.

The second statement immediately follows from the definition of
$\mathcal T$.
\end{proof}
\medskip

\begin{rem}
{\rm
Consider now  the projection map $\pi^\diff:\diff\to\W$ in the
presence of the distribution ${\mathcal T}$ on $\diff$. The path
lifting property in this case is a restatement of the nonholonomic
Moser theorem. Namely, for a curve
$\{\mu_t~|~\mu_0=\mu\}$ in the space $\W$ of smooth densities
Theorem  \ref{thm:nonhol_vol}
proves that there is a curve $\{g^t~|~g^0=id\}$ in $\diff$,
everywhere tangent to the distribution  ${\mathcal T}$ and
projecting to  $\{\mu_t\}$: $\pi^\diff(g^t)=\mu_t$.

Recall that in the classical case the corresponding path
lifting becomes unique once we fix the gradient horizontal bundle
$Hor_\phi\subset T_\phi\diff$ for any diffeomorphism $\phi\in\diff$.
Similarly, in the nonholonomic case we consider the spaces of gradient projections
instead of the gradient spaces:
$H\!or^\tau_{id}:=\{P^\tau \nabla f~|~f\in C^\infty(M)\}$,
where $P^\tau$ stands for the orthogonal projection onto the
distribution $\tau$ in a given Riemannian metric on $M$. The
right-translated gradient projections $H\!or^\tau_\phi:=\{(P^\tau
\nabla f)\circ\phi~|~f\in C^\infty(M)\}$ define a horizontal bundle for the principal
bundle $\diff\to\W$ by nonholonomic Hodge decomposition. (Note also
that in both classical and nonholonomic cases, the obtained
horizontal distributions on $\diff$ are nonintegrable, cf. \cite{Ot}.
Indeed, the Lie bracket of two gradient fields is not necessarily
a gradient field, and similarly for gradient projections. Hence there
are no horizontal sections of the bundle $\diff\to\W$, tangent to
these horizontal gradient distributions.)

\medskip

As we will see in Sections \ref{sec:geom} and \ref{sec:subriem},
both  gradient fields $\{\nabla
f\}$ in the classical case and gradient projections $ \{P^\tau
\nabla f\}$ in the nonholonomic case allow one to move the
densities in the ``fastest way'', and are important in transport
problems of finding optimal (``shortest") path between densities.
}
\end{rem}

\medskip

%%%%%%%%%%%%%%%%%%%%%%

\subsection{Accessibility of diffeomorphisms and symplectic structures}
Presumably, even a stronger statement holds:

\begin{conj}\label{conj:access}
Every diffeomorphism  in  the diffeomorphism group $\diff$ can be
accessed by a horizontal path tangent to the distribution $\mathcal
T$.
\end{conj}

This conjecture can be thought of as an analog of the Chow-Rashevsky
theorem in the infinite-dimensional setting of the group of
diffeomorphisms, provided that the distribution $\mathcal T$ is
bracket-generating on $\diff$. Note, however, that the
Chow-Rashevsky theorem  is unknown in the general setting of an
infinite-dimensional manifold, while there are only ``approximate"
analogs of it, e.g. on a Hilbert manifold.

A proof of this conjecture on accessibility of all diffeomorphisms
by flows of vector fields tangent to a nonholonomic distribution
would imply the nonholonomic Moser theorem \ref{thm:nonhol_vol} on
volume forms. Moreover, it would also imply the following nonholonomic
version of the Moser theorem on symplectic structures from
\cite{Moser1965}.

\begin{conj} Suppose that on a manifold $M$ two symplectic structures
$\omega_0$ and $\omega_1$ from the same cohomology class can be
connected by a path of symplectic structures in the same class.
Then for a  bracket-generating
distribution $\tau$ on $M$ there exists a diffeomorphism $\phi$ of
$M$ which is the time-one-map of a non-autonomous vector field $V_t$
tangent to the distribution $\tau$ everywhere on $M$ and for every
$t\in [0,1]$, such that $\phi^*\omega_1=\omega_0$.
\end{conj}

This conjecture follows from the one above since one would consider
the diffeomorphism from the classical Moser theorem, and realize it
by the horizontal path (tangent to the distribution $\mathcal T$) on
the diffeomorphism group, which  exists if Conjecture
\ref{conj:access} holds.

\bigskip

%%%%%%%%%%%%%%%%%%%%%%%%%%%%%%%%%%%%%%%%%%%%%%%%

%%%%%%%%%%%%%%%%%%%%%%%%%%%%%%%%%%%%%%%%%%%%%%%%

\section{The Riemannian geometry of diffeomorphism groups and mass transport}\label{sec:geom}

The differential geometry of diffeomorphism groups is closely
related to the theory of optimal mass transport, and in particular,
to the problem of moving one density to another while minimizing
certain cost on a Riemannian manifold. In this section,
we review the corresponding metric properties of the diffeomorphism
group and the space of volume forms.

\subsection{Optimal mass transport}
Let $M$ be a compact Riemannian manifold without boundary (or, more
generally, a complete metric space) with a distance function $d$.
Let $\mu$ and $\nu$ be two Borel probability measures on the
manifold $M$ which are absolutely continuous with respect to the
Lebesgue measure. Consider the following optimal mass transport
problem: Find a Borel map $\phi:M\to M$ that pushes the measure
$\mu$ forward to $\nu$ and attains the infimum of the $L^2$-cost
functional $\int_M d^2(x,\phi(x))\mu$ among all such maps.

The set of all Borel probability measures is called the Wasserstein
space. The minimal cost of transport defines a metric $\td$ on this
space:
\begin{equation}\label{optimal}
{\td}^2(\mu, \nu):=\inf_\phi\{\int_M d^2(x,\phi(x))\mu~|~\phi_*\mu=\nu\}\,.
\end{equation}

This mass transport problem admits a unique solution
$\phi$ (defined up to measure zero sets), called an optimal map (see
\cite{Br} for $M=\Real^n$ and \cite{Mc} for any compact connected
Riemannian manifold $M$ without boundary). Furthermore, there exists
a 1-parameter family of Borel maps $\phi_t$ starting at the identity
map $\phi_0=id$, ending at the optimal map $\phi_1=\phi$ and such that
$\phi_t$ is the optimal map pushing $\mu$ forward to $\nu_t:=\phi_{t*}\mu$
for any $t\in(0,1)$. The corresponding 1-parameter family of
measures $\nu_t$ describes a geodesic in the Wasserstein
space of measures with respect to the distance function $\td$ and is
called the displacement interpolation between $\mu$ and $\nu$, see
\cite{Vi1} for details. (More generally, in mass transport problems
one can replace $d^2$ in the above formula by a cost function
$c:M\times M\to \Real$, while we mostly focus on the case $c=d^2/2$
and its subriemannian analog.)

In what follows, we consider  a smooth version of the
Wasserstein space, cf. Section \ref{subsec:fibra}. Recall that the smooth
Wasserstein space $\W$ consists of smooth volume forms with the total
integral equal to 1. One can consider an infinite-dimensional manifold
structure on the smooth Wasserstein space, a (weak) Riemannian
metric $\langle \,,\,\rangle ^\W$, corresponding to the distance function
$\td$, and geodesics on this space. Similar to the finite-dimensional case,
geodesics on the smooth Wasserstein space $\W$ can be formally defined as
projections of trajectories of the Hamiltonian vector field with
the ``kinetic energy'' Hamiltonian in the tangent bundle $T\W$.

\bigskip
%%%%%%%%%%%%%%%%%%%%%%%%%%%%%%%

\subsection{The Otto calculus}
For a Riemannian manifold $M$ both spaces $\mathcal D$ and $\W$ can
be equipped with (weak) Riemannian structures, i.e. can be formally
regarded as infinite-dimensional Riemannian manifolds, cf.
\cite{EbMa}. (One can consider $H^s$-diffeomorphisms and
$H^{s-1}$-forms of Sobolev class $s> n/2 +1$. Both sets can be
considered as smooth Hilbert manifolds. However, this is not
applicable in the subriemannian case, discussed later, hence we confine
to the $C^\infty$ setting applicable in the both cases.)

From now on we fix a Riemannian metric $\langle ,\rangle ^M$ on the
manifold $M$, whose Riemannian volume is the form $\mu$.
On the diffeomorphism group we define a Riemannian metric $\langle ,\rangle ^\diff$ whose
value at a point $\phi\in \diff$ is given by
\begin{equation}\label{dmetric}
\langle X_1\circ\phi,X_2\circ\phi\rangle ^\diff
:=\int_M\langle X_1\circ\phi(x),X_2\circ\phi(x)\rangle ^M_{\phi(x)}\mu.
\end{equation}
The action along a curve (or, ``energy'' of a curve)
$\{\phi_t~|~t\in[0,1]\}\subset \diff$ in this metric is defined
in the following straightforward way:
$$
E(\{\phi_t\})=\int_0^1dt\int_M \langle \partial_t\phi_t, \partial_t\phi_t\rangle ^M
\,\mu\, .
$$
%(For $H^s$ maps $\phi$ with $s>n/2+1$ this metric endows $\diff$ with
%the structure of a Hilbert manifold, see details in \cite{EbMa}.)
If $M$ is flat, $\diff$ is locally isometric to the (pre-)Hilbert $L^2$-space
of (smooth) vector-functions $\phi$, see e.g. \cite{Shnir}.
The following proposition is well-known.

%If $M$ is flat, this metric locally identifies $\diff$ with the
%Hilbert $L^2$-space of $C^s$-functions $\phi(x)$ with $s>n/2+1$ in
%a vicinity of the identity map $\phi_0(x)=x$, see details in
%\cite{EbMa}. Let $\nabla$ be the Levi-Civita connection on the manifold $M$.

\begin{prop}\label{prop:burg}
Let $\phi_t$ be a geodesic on the diffeomorphism group $\diff$ with
respect to the above Riemannian metric $\langle ,\rangle ^\diff$, and $V_t$ be the
(time-dependent) velocity field of the corresponding flow:
$\partial_t\phi_t=V_t\circ\phi_t$. Then the velocity $V_t$ satisfies
the inviscid Burgers equation on $M$:
\[
\partial_t V_t+\nabla_{V_t}V_t=0\,,
\]
where $\nabla_{V_t}V_t$ stands for the covariant derivative of the field $V_t$
on $M$ along itself.
\end{prop}

\begin{proof}
In the flat case the geodesic equation is $\partial_t^2{\phi_t}=0$:
this is the Euler-Lagrange equation for the action functional $E$.
Differentiate $\partial_t\phi_t=V_t\circ\phi_t$ with respect to time $t$
and use this geodesic equation to obtain
\begin{equation}\label{Burger1}
\partial_t V_t\circ\phi_t+\nabla_{V_t}\partial_t\phi_t=0.
\end{equation}
After another substitution  $\partial_t\phi_t=V_t\circ\phi_t$, the
later becomes
\[
(\partial_t V_t+\nabla_{V_t}V_t)\circ\phi_t=0,
\]
which is equivalent to the Burgers equation.

The non-flat case involves differentiation in the Levi-Civita connection
on $M$ and leads to the same Burgers equation, see details in \cite{EbMa, KM07}.
\end{proof}

\begin{rem}\label{rem:particles}
{\rm
Smooth solutions of the Burgers equation correspond to
non-interacting particles on the manifold $M$ flying along those geodesics
on $M$ which are defined by the initial velocities $V_0(x)$.
The Burgers flows have the form $\phi_t(x)=\exp^M(tV_0(x))$, where $\exp^M:TM\to M$
is the Riemannian exponential map on $M$.
}
\end{rem}

\medskip

\begin{prop}\cite{Ot}\label{rsub}
The bundle projection $\pi^\diff:\diff\to\W$ is a
Riemannian submersion of the metric $\langle \,,\,\rangle ^\diff$ on the diffeomorphism group
$\diff$ to the Riemannian metric $\langle \,,\,\rangle ^\W$ on the smooth Wasserstein space $\W$
for the $L^2$-cost.
The horizontal (i.e. normal to fibers) spaces in the bundle $\diff \to\W$
are right-translated gradient fields.
\end{prop}

Recall that for two Riemannian manifolds $Q$ and $B$, a Riemannian
submersion $\pi: Q\to B$ is a mapping onto $B$ which has maximal
rank and preserves lengths of horizontal tangent vectors to $Q$, see
e.g. \cite{Neill}. For a bundle $Q\to B$, this means that there is a
distribution of horizontal spaces on $Q$, orthogonal to the fibers,
which is projected isometrically to the tangent spaces to $B$. One of the
main properties of a Riemannian submersion gives the following feature
of geodesics:

\begin{cor}\label{horizontal}
Any geodesic, initially tangent to a horizontal space on the full
diffeomorphism group $\diff$, always remains horizontal, i.e.  tangent to
the horizontal distribution.  There is a one-to-one correspondence between
geodesics on the base $\W$ starting at the measure $\mu$ and horizontal
geodesics in $\diff$ starting at the identity diffeomorphism $id$.
\end{cor}

\begin{rem}\label{MA1}
{\rm
In the PDE terms, the horizontality of a geodesic means that a
solution of the Burgers equation with a potential initial condition
remains potential forever. This also follows from the
Hamiltonian formalism and the moment map geometry discussed in the
next section.
Since horizontal geodesics in the group $\diff$ correspond to geodesics
on the density space $\W$, potential
solutions of the Burgers equation (corresponding to horizontal
geodesics) move the densities in the fastest way. The corresponding
time-one-maps for Burgers potential solutions provide optimal maps
for moving the density $\mu$ to any other density $\nu$, see \cite{Br,
Mc}.

The Burgers potential solutions have the form $\phi_t(x)=\exp^M(-t\nabla
f(x))$ as long as the right-hand-side is smooth. The time-one-map
$\phi_1$ for the flow $\phi_t$ provides an optimal map
between probability measures if the function $f$ is a $(d^2/2)$-{\it concave}
function. The notion of $c$-concavity for a cost function $c$ on $M$ is defined as
follows. For a function $f$ its $c$-transform is
$f^c(y)=\inf_{x\in M}(c(x,y)-f(x))$ and the function $f$ is said to be
$c$-concave if $f^{cc}=f$. Here, we consider the case $c=d^2/2$.
The family of maps $\phi_t$ defines the displacement interpolation
mentioned in Section 4.1.

Let $\theta$ and $\nu$ be volume forms with the same total volume
and let $g$ and $h$ be functions on the manifold $M$ defined by $\theta=g\vol$ and
$\nu=h\vol$, where $\vol$ be the Riemannian volume form.
Then a diffeomorphism $\phi$ moving one density to the other ($\phi_*\theta=\nu$)
satisfies $h(\phi(x))\det(D\phi(x))=g(x)$, where $D\phi$ is the Jacobi matrix of
the diffeomorphism $\phi$. In the flat case
the optimal map $\phi$ is gradient, $\phi=\nabla \tilde f$,
and the corresponding convex potential $\tilde f$
satisfies the Monge-Amp\`{e}re equation
\[
\det(\operatorname{Hess}\tilde f(x)))=\frac{g(x)}{h(\nabla \tilde f(x))}\,,
\]
since $D(\nabla \tilde f)=\operatorname{Hess}\tilde f$.
In the non-flat case, the optimal map is $\phi(x)=\exp^M(-\nabla f(x))$
for a $(d^2/2)$-concave potential $ f$, and the equation is Monge-Amp\`{e}re-like,
see \cite{Mc, Vi1} for details. Below we describe the corresponding nonholonomic
analogs of these objects.
}
\end{rem}

\bigskip

%%%%%%%%%%%%%%%%%%%%%%%%%%%%%%%

%%%%%%%%%%%%%%%%%%%%%%%%%%%%%%%

\section{The Hamiltonian mechanics on  diffeomorphism groups}

In this section we present a Hamiltonian framework for the Otto
calculus and, in particular, give a symplectic proof of Proposition \ref{rsub} and
Corollary \ref{horizontal} on the submersion properties
along with their generalizations.

\subsection{Averaged Hamiltonians}
We fix a Riemannian metric $\langle ,\rangle ^M$ on the manifold $M$ and consider
the corresponding Riemannian metric  $\langle ,\rangle ^\diff$ on the diffeomorphism
group $\diff$. This defines a map $(X\circ\phi)\mapsto \langle X\circ\phi\,,\cdot\rangle ^\diff$
from the tangent bundle $T\diff$ to the cotangent bundle $T^*\diff$. By using
this map, one can pull back the canonical symplectic form $\omega^{T^*\diff}$
from the cotangent bundle $T^*\diff$ to the tangent bundle $T\diff$, and
regard the latter as a manifold equipped with the symplectic form
$\omega^{T\diff}$.\footnote{The consideration of the tangent bundle
$T\diff$ (instead of $T^*\diff$) as a symplectic manifold allows one
to avoid dealing with duals of infinite-dimensional spaces here.}
Similarly, a symplectic structure $\omega^{TM}$ can be defined on the tangent
bundle $TM$ by pulling back the canonical symplectic form on the cotangent
bundle $T^*M$ via the Riemannian metric $\langle ,\rangle ^M$.
The two symplectic forms are related as follows. A tangent vector
$V$ in the tangent space $T_{X\circ\phi}T\diff$
at the point $X\circ\phi\in T\diff$ is a map from $M$ to $T(TM)=T^2M$
such that $\pi^{T^2M}\circ V=X\circ\phi$, where $\pi^{T^2M}:T(TM)\to TM$
is the tangent bundle projection. Let $V_1$ and $V_2$ be two tangent vectors
in  $T_{X\circ\phi}T\diff$ at the point $X\circ\phi$, then the symplectic forms
are related in the following way:
$$
\omega^{T\diff}(V_1,V_2)=\int_M\omega^{TM}(V_1(x),V_2(x))\mu(x)\,,
$$
where $\omega^{TM}$ is understood as the pairing on  $T(TM)=T^2M$.

\begin{defn}
{\rm
Let $H^M$ be a Hamiltonian function on the tangent bundle $TM$ of
the manifold $M$. The {\it averaged Hamiltonian function} is the
function $H^\diff$ on the tangent bundle $T\diff$ of the
diffeomorphism group $\diff$ obtained by averaging the corresponding
Hamiltonian $H^M$ over $M$ in the following way: its value at a
point $X\circ\phi\in T_\phi\diff$ is
\begin{equation}\label{eq:Ham}
H^\diff(X\circ\phi):=\int_MH^M(X\circ\phi(x))\mu(x)\,
\end{equation}
for a vector field $X\in\vect$ and a diffeomorphism $\phi\in \diff$.
}
\end{defn}

Consider the Hamiltonian flows for these Hamiltonian functions $H^M$
and $H^\diff$ on the tangent bundles $TM$ and  $T\diff$,
respectively, with respect to the standard symplectic structures on
the bundles. The following theorem can be viewed as a generalization
of Propositions \ref{prop:burg} and \ref{rsub}.

\begin{thm}\label{genham}
Each Hamiltonian trajectory for the averaged Hamiltonian function
$H^\diff$ on $T\diff$ describes a flow on the tangent bundle $TM$,
in which every tangent vector to $M$ moves along its own
$H^M$-Hamiltonian trajectory in $TM$.
\end{thm}

\begin{exmp}\label{ex:kinetic}
{\rm
For the Hamiltonian $K^M(p,q)=\frac{1}{2}\langle p,p\rangle ^M$ given by the
``kinetic energy" for the metric on $M$, the above theorem implies
that any geodesic on $\diff$ is a family of diffeomorphisms of $M$,
in which each particle moves along its own geodesic on $M$ with
constant velocity, i.e. its velocity field is a solution to the
Burgers equation, cf. Remark \ref{rem:particles}.
}
\end{exmp}

Below we discuss this theorem and its geometric meaning in detail.
In particular, in the above form, the statement is  also applicable
to the case of nonholonomic distributions (i.e. subriemannian, or
Carnot-Caratheodory spaces) discussed in the next section.
\medskip

%%%%%%%%%%%%%%%%%%%%

\subsection{Riemannian submersion and symplectic quotients}
\medskip

We start with a Hamiltonian proof of Proposition \ref{rsub}
on the Riemannian submersion $\diff\to\W$ of diffeomorphisms onto densities.
Recall the following general construction in symplectic geometry.
Let $\pi:Q\to B$ be a principal bundle with the structure group $G$.

\begin{lem} (see e.g. \cite{Cannas})
The symplectic reduction of the cotangent bundle $T^*Q$ over the $G$-action
gives the cotangent bundle $T^*B=T^*Q/\!/G$.
\end{lem}

\begin{proof}
The moment map $J:T^*Q\to\Lg^*$ associated with this action takes
$T^*Q$ to the dual of the Lie algebra $\Lg=Lie(G)$. For the
$G$-action on   $T^*Q$ the moment map $J$ is the projection of any
cotangent space $T_a^*Q$ to cotangent space $T_a^*F\approx\Lg^*$
for the fiber $F$ through a point $a\in Q$. The preimage
$J^{-1}(0)$ of the zero value is the subbundle of $T^*Q$
consisting of covectors vanishing on  fibers.
Such covectors are naturally identified with covectors on the base $B$.
Thus factoring out the $G$-action, which moves the point $a$ over the fiber $F$, 
we obtain the bundle $T^*B$.
\end{proof}

Suppose also that $Q$ is equipped with a $G$-invariant Riemannian metric $\langle ,\rangle ^Q$.

\begin{lem}
The Riemannian submersion of $(Q, \langle ,\rangle ^Q)$ to the base $B$ with the
induced metric $\langle ,\rangle ^B$ is the result of the symplectic reduction.
\end{lem}

\begin{proof}
Indeed, the metric $\langle ,\rangle ^Q$ gives a natural identification
$T^*Q\approx TQ$ of the tangent and cotangent bundles for $Q$, and
the ``projected  metric'' is equivalent to a similar identification for
the base manifold $B$.

In the presence of metric in $Q$, the preimage $J^{-1}(0)$ is identified with
all vectors in $TQ$ orthogonal to fibers, that is $J^{-1}(0)$ is the
horizontal subbundle in $TQ$. Hence, the symplectic quotient $J^{-1}(0)/G$
can be identified with the tangent bundle $TB$.
\end{proof}

\medskip

\proofbegin of Proposition \ref{rsub}.
Now we apply this ``dictionary" to the diffeomorphism group $\diff$ and the
Wasserstein space $\W$. Consider the projection map
$\pi^\diff:\diff\to\W$ as a principal bundle with the structure
group $\sdiff$ of volume-preserving diffeomorphisms of $M$. Recall that the
vertical space of this principal bundle at a point $\phi\in \diff$ consists of
right-translations by the diffeomorphism $\phi$ of vector fields which are
divergence-free with respect to the volume form $\phi_*\mu$:
$
V\!er_\phi=\{X\circ\phi~|~\di_{\phi_*\mu}X=0\}\,,
$ and the horizontal space is given by
translated gradient fields:
$
H\!or_\phi=\{\nabla f\circ\phi~|~f\in C^\infty(M)\}.
$

For each volume-preserving diffeomorphism $\psi\in\sdiff$,
the $\sdiff$-action $R_\psi$ of $\psi$ by
right translations on the diffeomorphism group is given by
\[
R_\psi(\phi)=\phi\circ\psi.
\]
The induced action
$TR_\psi:T\diff\to T\diff$ on the tangent spaces of the
diffeomorphism group is given by
\[
TR_\psi(X\circ\phi)=(X\circ\phi)\circ\psi.
\]
One can see that for volume-preserving diffeomorphisms $\psi$ this
action preserves the Riemannian metric (\ref{dmetric}) on the
diffeomorphism group $\diff$ (it is the change of variable
formula), while for a general diffeomorphism one has an extra
factor $D\psi$, the Jacobian of $\psi$, in the integral.
\proofend

\bigskip

\begin{rem}
{\rm
The explicit formula of the moment map $J:TQ\to\svect^*$
for the group of volume-preserving
diffeomorphisms $G=\sdiff$ acting on $Q=\diff$ is
\[
J(X\circ\phi)(Y)=\int_M\langle X,\phi_*Y\rangle ^M\phi_*\mu\,,
\]
where $Y\in \svect$ is any
vector field  on $M$ divergence-free with respect to the volume form $\mu$,
$X\in\vect$,  and $\phi\in\diff$.
}
\end{rem}

\bigskip
%%%%%%%%%%%%%%%%%%%%%%%%%%%%

\subsection{Hamiltonian flows on the diffeomorphism groups}
%{Geodesics and the reduction}
Let $H^Q:TQ\to\Real$ be a Hamiltonian function invariant under
the $G$-action on the cotangent bundle of the total space $Q$.
The restriction of the function $H^Q$ to the horizontal bundle
$J^{-1}(0)\subset TQ$ is also $G$-invariant, and hence
descends to a function $H^B:TB\to\Real$  on the symplectic quotient, the
tangent bundle of the base  $B$. Symplectic quotients admit
the following reduction of Hamiltonian dynamics:

\begin{prop}\cite{Cannas}\label{redy}
The Hamiltonian flow of the function $H^Q$ preserves the preimage
$J^{-1}(0)$, i.e. trajectories with horizontal initial conditions
stay horizontal. Furthermore, the Hamiltonian flow of the function
$H^Q$ on the tangent bundle $TQ$ of the total space $Q$ descends
to the Hamiltonian flow of the function $H^B$ on the
tangent bundle $TB$ of the base.
\end{prop}

Now we are going to apply this scheme to the bundle $\diff\to\W$.
For  a fixed Hamiltonian function $H^M$ on the tangent bundle $TM$
to the manifold $M$, consider the corresponding averaged Hamiltonian
function $H^\diff$ on $T\diff$, given by the formula (\ref{eq:Ham}):
$H^\diff(X\circ\phi):=\int_MH^M(X\circ\phi(x))\mu$. The latter
Hamiltonian is $\sdiff$-invariant (as also follows from the change
of variable formula) and it will play the role of the function
$H^Q$. Thus the flow for the averaged Hamiltonian $H^\diff$ descends to the flow
of a certain Hamiltonian $H^\W$ on $T\W$.

Describe explicitly the corresponding flow on
the tangent bundles of $\diff$ and $\W$. Let $\Psi^{H^M}_t:TM\to
TM$ be the Hamiltonian flow of the Hamiltonian $H^M$ on the tangent
bundle of the manifold $M$ and $\Psi^{H^\diff}_t:T\diff\to T\diff$
denotes the flow for the Hamiltonian function $H^\diff$ on the
tangent bundle of the diffeomorphism group.

\begin{thm}\label{HamFl}{\bf (=\ref{genham}$'$)}
The Hamiltonian flows of the Hamiltonians $H^\diff$ and $H^M$ are
related by
\[
\Psi^{H^\diff}_t(X\circ\phi)(x)=\Psi^{H^M}_t( X(\phi(x)))\,,
\]
where, on the right-hand-side, the flow $\Psi^{H^M}_t$ on $TM$ transports the shifted
field $X(\phi(x))$, while, on the
left-hand-side, $X\circ\phi$ is regarded as a tangent vector to
$\diff$ at the point $\phi$.
\end{thm}

\begin{proof}
Prove this infinitesimally (cf. \cite{EbMa}).
Let $X_{H^\diff}$ and $X_{H^M}$ be the Hamiltonian
vector fields corresponding to the Hamiltonians $H^\diff$ and $H^M$ respectively.
We claim that $X_{H^\diff}(X\circ\phi)=X_{H^M}\circ X\circ\phi$.
Indeed, by the definition of Hamiltonian fields, we have
$$
\omega^{T\diff}(X_{H^M}\circ X\circ\phi,Y)
=\int_M\omega^{TM}(X_{H^M}(X(\phi(x))),Y(x))\mu
=\int_M dH^M_{X(\phi(x))}(Y(x))\mu(x)\,
$$
for any $Y\in T_\phi\diff$.
By interchanging the integration and exterior differentiation, the latter
expression becomes $dH^\diff_{X\circ\phi}(Y)$ and the result follows.
\end{proof}

\begin{rem}
{\rm
This theorem has a simple geometric meaning for
the ``kinetic energy'' Hamiltonian function $K^M(v):=\frac 12\langle v,v\rangle ^M$ on the
tangent bundle $TM$. One of the possible definitions of geodesics in $M$
is that they are projections to $M$ of trajectories of the
Hamiltonian flow on $TM$, whose Hamiltonian function is the kinetic
energy. In other words, the Riemannian exponential map $\exp^M$ on
the manifold $M$ is the projection of the Hamiltonian flow
$\Psi^{K^M}_t$ on $TM$. Similarly, the Riemannian exponential
$\exp^\diff$ of the diffeomorphism group $\diff$ is the projection
of the Hamiltonian flow for the Hamiltonian
$K^\diff(X\circ\phi):=\frac 12\int_M\langle X\circ\phi,X\circ\phi\rangle^M \mu$ 
on $T\diff$.
%\[\exp^M(X(x))=\pi^{TM}(\Psi^{H^M}_t(X(x))).\]

%A simple calculation shows that the two exponential maps are related by
%\[\exp^\diff(X)=\exp^M\circ X.\]

Recall that the geodesics on the diffeomorphism group (described by
the Burgers equation, see Proposition \ref{prop:burg}) starting at
the identity with the initial velocity $V\in T_{id}\diff$ are the
flows which move each particle $x$ on the manifold $M$ along the
geodesic with the direction $V(x)$. Such a geodesic is well defined
on the diffeomorphism group $\diff$ as long as the particles do not
collide. The corresponding Hamiltonian flow on the tangent bundle
$T\diff$ of the diffeomorphism group describes how the corresponding
velocities of these particles vary (cf. Example \ref{ex:kinetic}).
\medskip

For a more general Hamiltonian $H^M$ on the tangent bundle $TM$, each
particle $x\in M$ with an initial velocity $V(x)$ will be moving along
the corresponding characteristic, %$\pi^{TM}(\Psi_t^{H^M}(X(x)))$,
which is the projection to $M$ of the corresponding trajectory
$\Psi_t^{H^M}(V(x))$ in the tangent bundle $TM$.
%(Here $\pi^{TM}:TM\to M$ is the bundle projection.)
}
\end{rem}

Now we would like to describe more explicitly horizontal geodesics
and characteristics on the diffeomorphism group
$\diff$. Recall that $\Psi^{H^\diff}_t$ denotes the Hamiltonian flow
of the averaged Hamiltonian $H^\diff$ on the tangent bundle $T\diff$
of the diffeomorphism group $\diff$.
If this Hamiltonian flow is gradient at the initial moment, it always stays gradient,
as implied by Corollary \ref{horizontal}.
Furthermore, the corresponding potential can be described as follows.

\begin{cor}\label{HJ}
Let $f$ be a function on the manifold $M$. Then the Hamiltonian flow
for $H^\diff$ with the initial condition $\nabla f\circ\phi\in
T_\phi\diff$ has the form $\nabla f_t\circ\phi_t$, where $\phi_t\in\diff$
is a family of diffeomorphisms and $f_t$ is
the family of functions on $M$ starting at $f_0=f$ and satisfying
the Hamilton-Jacobi equation
\begin{equation}\label{HJeqn}
\partial_t f_t+H^M(\nabla f_t(x))=0\,.
\end{equation}
\end{cor}

\begin{figure}[ht!]
\input epsf
\centerline{\epsfysize=0.30\vsize\epsffile{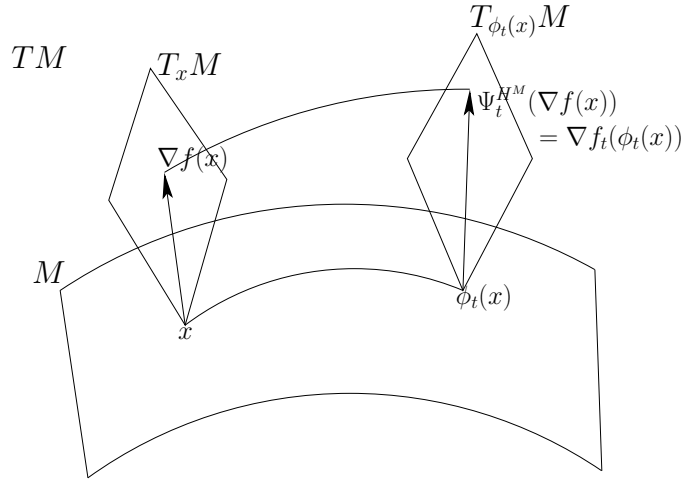}}
 \caption{Hamiltonian flow of the Hamiltonian $H^M$ and its projection: 
The curve $\phi_t(x)$ is the projection of the curve 
$\Psi_t^{H^M}(\nabla f(x))$ to the manifold $M$.}\label{fig:projHam}
\end{figure}

\begin{proof}
This follows from the method of characteristics, which gives
the following way of finding
$f_t$, the solution to the Hamilton-Jacobi equation (\ref{HJeqn}).
Consider the tangent vector $\nabla f(x)$ for each point $x\in M$.
Denote by $\Psi_t^{H^M}:TM\to TM$ the Hamiltonian flow for the Hamiltonian
$H^M:TM\to\Real$ and consider its trajectory $t\mapsto\Psi_t^{H^M}(\nabla f(x))$
starting at the tangent vector $\nabla f(x)$.
Then project this trajectory to $M$ using the tangent bundle projection
$\pi^{TM}:TM\to M$ to obtain a curve in $M$.
It is given by the formula $t\mapsto\pi^{TM}(\Psi_t^{H^M}(\nabla f(x)))$.
As $x$ varies over the manifold $M$, this defines a flow
$\phi_t:=\pi^{TM}\circ\Psi_t^{H^M}\circ\nabla f$ on $M$.
(Note that this procedure defines a flow for small time $t$, while
for larger times the map $\phi_t$ may cease to be a diffeomorphism, 
i.e. shock waves can appear.)
The corresponding time-dependent vector field is gradient and defines the family
$\nabla f_t$, the gradient of the solution to the Hamilton-Jacobi equation above,
see Figure \ref{fig:projHam}.
\end{proof}

%\begin{figure}[ht!]
%\begin{center}
%\centerline{\scalebox{0.5}{\input{mostranfig4.pstex_t}}}
% \caption{Hamiltonian flow of the Hamiltonian $H^M$ and its projection: 
%The curve $\phi_t(x)$ is the projection of the curve 
%$\Psi_t^{H^M}(\nabla f(x))$ to the manifold $M$.}\label{fig:projHam}
%\end{center}
%\end{figure}

\begin{rem}\label{rem:charact}
{\rm
The above corollary manifests that the Hamilton-Jacobi
equation (\ref{HJeqn}) can
be solved using the method of characteristics due to the built-in symmetry
group of all volume preserving diffeomorphisms.
}
\end{rem}
%%%%%%%%%%%%%%%%%%%%%%%%%%%%%%%%%%%%%%%%%%%%%%%%%%%

\subsection{Hamiltonian flows on the Wasserstein space}
What is the corresponding flow on the tangent bundle $T\W$ of the
Wasserstein space, induced by the Hamiltonian flow on $T\diff$ for
the diffeomorphism group $\diff$ after the projection $\pi^\diff:\diff\to\W$?
Fix a Hamiltonian $H^M$ on the tangent bundle $TM$ which defines the
averaged Hamiltonian function $H^\diff$ on the tangent bundle $T\diff$, see
Equation (\ref{eq:Ham}). Describe explicitly the induced Hamiltonian $H^\W$
on the tangent bundle $T\W$.

Let $(\nu,\eta)$ be a tangent vector at a density
$\nu$ on $M$, regarded as a point of the Wasserstein space $\W$.
The normalization of densities ($\int\nu=1$ for all $\nu\in \W$)
gives the constraint for tangent vectors: $\int_M\eta=0$.
Let $f:M\to\Real$ be a function that satisfies $(-\di_\nu{\nabla f})\nu=\eta$.
(Given  $(\nu,\eta)$, such a function is defined uniquely up to 
an additive constant.)
Then the induced Hamiltonian on the tangent bundle $T\W$ of the base $\W$ is
given by
\begin{equation}\label{reHam}
H^\W(\nu,\eta)=\int_M H^M(\nabla f(x))\,\nu\,,
\end{equation}
since $\nabla f$ is a vector of the horizontal distribution in $T\diff$.

Now, the flow $\Psi_t^{H^\W}$ of the corresponding Hamiltonian field
on $T\W$ can be found explicitly by employing  Proposition \ref{redy}.
Consider the flow $\phi_t:=\pi^{TM}\circ\Psi_t^{H^M}\circ\nabla f$
defined on $M$ for small $t$ in Corollary \ref{HJ}.

\begin{thm}\label{main}
The Hamiltonian flow $\Psi_t^{H^\W}$ of the Hamiltonian function
$H^\W$ on the tangent bundle $T\W$ of the Wasserstein space $\W$
is
\[
\Psi_t^{H^\W}(\nu,\eta)=(\nu_t,-{\LD}_{\nabla f_t}\nu_t)\,,
\]
where $\LD$ is the Lie derivative, the family of functions $f_t$ satisfies
the Hamilton-Jacobi equation (\ref{HJeqn}) for the
Hamiltonian function $H^M$ on the tangent bundle $TM$, and the
family $\nu_t=(\phi_t)_*\nu$ is the push forward of the volume
form $\nu$ by the map $\phi_t$ defined above.
\end{thm}

\begin{proof}
The function $H^{\diff}(X\circ\phi)=\int_MH^M(X(\phi(x)))\mu(x)$ on the tangent
bundle $T\diff$ of the diffeomorphism group induces the Hamiltonian $H^\W$ on $T\W$.
By virtue of the Hamiltonian reduction,  Hamiltonian trajectories of $H^\diff$
contained in the horizontal bundle $H\!or=\{\nabla f\circ\phi~|~f\in C^\infty(M)\}$
descend to Hamiltonian trajectories of $H^\W$. Then the Hamiltonian flow
$\Psi^{H^\diff}$ of the Hamiltonian $H^\diff$ is given by
$\Psi^{H^\diff}(X\circ\phi)=\Psi^{H^M}\circ X\circ\phi$, due to Theorem \ref{HamFl}.
By restricting this to the horizontal bundle $H\!or$ we have
\begin{equation}\label{HamFlHor}
\Psi^{H^\diff}(\nabla f\circ\phi)=\Psi^{H^M}\circ \nabla f\circ\phi.
\end{equation}
The flow $\Psi^{H^\diff}$ is described in Corollary \ref{HJ} and has the form
$\Psi^{H^\diff}(\nabla f\circ\phi)=\nabla f_t\circ\phi_t$,
where $f_t$ and $\phi_t$ are defined as required.

On the other hand, recall that the projection $\pi^\diff:\diff\to\W$ is defined by
$\pi^\diff(\phi)=\phi_*\mu$. The differential $D\pi^\diff$ of this map $\pi^\diff$
is
$$
D\pi(X\circ\phi):=(\phi_*\mu,-\LD_X(\phi_*\mu))\,.
$$
The application of this relation to (\ref{HamFlHor}) gives the result.
\end{proof}

\begin{rem}\label{MA2}
{\rm
The time-one-map  for the above density flow $\nu_t$ in the
Wasserstein space $\W$ formally describes optimal transport maps for
the Hamiltonian $H^M$. In particular, it recovers
the optimal map recently obtained in \cite{BeBu}. One considers
the optimal transport problem for the functional
\[
\inf_\phi\{\int_M c(x,\phi(x))\mu~|~\phi_*\mu=\nu\}
\]
with the cost function $c$ defined by
$$
c(x,y)=\inf_{\{\gamma
\text{ paths between $x$ and
$y$}\}}\int_0^1 L(\gamma,\dot\gamma)\,dt\,,
$$
where the infimum is taken over paths $\gamma$ joining  $x$ and
$y$ and the Lagrangian $L:TM\to\Real$ satisfies certain regularity and convexity
assumptions, see \cite{BeBu}.
The corresponding Hamiltonian $H^M$ in Theorem \ref{main} is the 
Legendre transform of the Lagrangian $L$.
Note that for the ``kinetic energy'' Lagrangian  $K^M$,
the above map becomes the optimal map $\exp^M(-\nabla f)$
mentioned at the beginning of this section, with $\exp^M:TM\to M$
being the Riemannian exponential of the manifold $M$.
}
\end{rem}

\bigskip

%%%%%%%%%%%%%%%%%%%%%%%%%%%%%%%%%

%%%%%%%%%%%%%%%%%%%%%%%%%%%%%%%%%

\section{The subriemannian geometry of  diffeomorphism groups}\label{sec:subriem}

In this section we develop the subriemannian setting for the
diffeomorphism group. In particular, we derive the geodesic equations for the
``nonholonomic Wasserstein metric," and describe nonholonomic
versions of the Monge-Amp\`ere  and heat equations.

\smallskip

Let $M$ be a manifold with a fixed distribution $\tau$ on it. 
Recall that a subriemannian metric is a positive definite inner product 
$\langle\,,\,\rangle^\tau$ on each plane of the distribution $\tau$ 
smoothly depending on a point in $M$. Such a metric can be defined by the 
bundle map $\I:T^*M\to\tau$, sending a covector $\alpha_x\in T_x^*M$
to the vector $V_x$ in the plane $\tau_x$ such that 
$\alpha_x(U)=\langle V_x,U\rangle^\tau$ on vectors $U\in \tau_x$. 
The subriemannian Hamiltonian $H^\tau:T^*M\to\Real$ is 
the corresponding fiberwise quadratic form:
\begin{equation}\label{subHam}
H^\tau(\alpha_x)=\frac{1}{2}\langle
V_x,V_x\rangle ^\tau.
\end{equation}
Let $\Psi_t^{H^\tau}$ be the Hamiltonian flow for time $t$ of the subriemannian 
Hamiltonian $H^\tau$ on $T^*M$, while $\pi^{T^*M}:T^*M\to M$ is the cotangent 
bundle projection. Then the subriemannian exponential map 
$\exp^\tau:T^*M\to M$ is defined as the projection to $M$ of the time-one-map of 
the above Hamiltonian flow on $T^*M$:
\begin{equation}\label{subexp}
\exp^\tau(t\alpha_x):=\pi^{T^*M}\Psi_t^{H^\tau}(\alpha_x).
\end{equation}
This relation defines a normal subriemannian geodesic on $M$ with
the initial covector $\alpha_x$. Note that the initial
velocity of the subriemannian geodesic $\exp^\tau(t\alpha_x)$ is
$V_x=\I \alpha_x\in\tau_x$. So, unlike the Riemannian case, there are many
subriemannian geodesics having the same initial velocity $V_x$ on $M$.

Let $d_\tau$ be a subriemannian (or, Carnot-Caratheodory) distance 
on the manifold $M$, defined as the infimum of the length of 
all absolutely continuous admissible (i.e. tangent to $\tau$)
curves joining given two points.
For a bracket-generating distribution $\tau$ any two points can be joined by
such a curve, so this distance is always finite. Consider the corresponding optimal
transport problem by replacing the Riemannian distance $d$
in (\ref{optimal}) with the subriemannian distance $d_\tau$.
Below we study the infinite-dimensional geometry of this subriemannian version
of the optimal transport problem.
Although in general normal subriemannian
geodesics might not exhaust all the length minimizing geodesics in
subriemannian manifolds (see \cite{Mo}), we will see that in the problems of
subriemannian optimal transport one can confine oneself to only such geodesics!

\medskip

\subsection{Subriemannian submersion}

Consider the following general setting: Let $(Q,\mathcal T)$ be a 
subriemannian space, i.e. a manifold $Q$ with a distribution $\mathcal T$ and 
a subriemannian metric $\langle\,,\,\rangle^\tau$ on it.
Suppose that $Q\to B$ is a bundle projection to 
a Riemannian base manifold $B$.

\begin{figure}[ht!]
\input epsf
\centerline{\epsfysize=0.4\vsize\epsffile{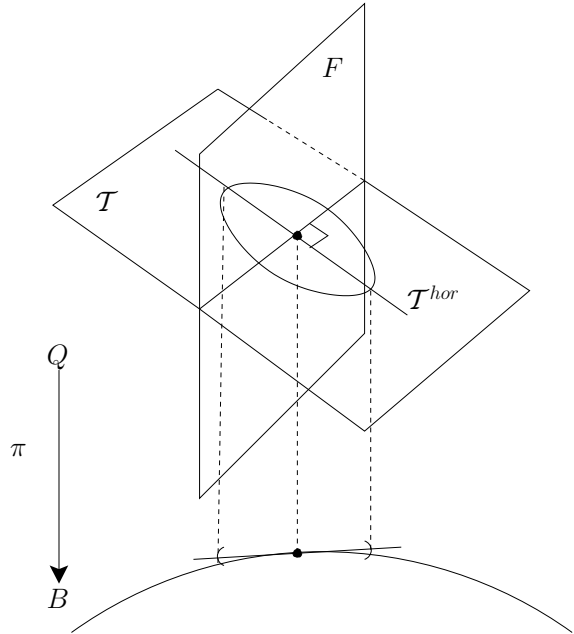}}
 \caption{Subriemannian submersion: horizontal subdistribution  $\mathcal T^{hor}$
is mapped isometrically to the tangent bundle $TB$ of the base.}\label{fig:subsub}
\end{figure}

\begin{defn}
{\rm
The projection $\pi: (Q,\mathcal T)\to B$ is a {\it subriemannian submersion}
if  the distribution $\mathcal T$ contains a {\it horizontal subdistribution}
 $\mathcal T^{hor}$, orthogonal (with respect to the subriemannian metric)
to the intersections of $\mathcal T$ with fibers, and the projection $\pi$
maps the spaces   $\mathcal T^{hor}$ isometrically to the tangent spaces of
the base $B$, see Figure \ref{fig:subsub}.} 
\end{defn}

Let a subriemannian submersion $\pi: (Q,\mathcal T)\to B$ be
a principal $G$-bundle  $Q\to B$,
where the distribution $\T$ and the subriemannian metric are invariant 
with respect to the action of the group $G$.
The following theorem is an analog of Corollary \ref{horizontal}.

%\begin{figure}[ht!]
%\begin{center}
%\centerline{\scalebox{0.5}{\input{mostranfig3.pstex_t}}}
% \caption{Subriemannian submersion: horizontal subdistribution  $\mathcal T^{hor}$
%is mapped isometrically to the tangent bundle $TB$ of the base.}\label{fig:subsub}
%\end{center}
%\end{figure}

\bigskip
\begin{thm}\label{thm:sub-reduction}
For each point $b$ in the base $B$ and a point $q$ 
in the fibre $\pi^{-1}(b)\subset Q$ over $b$, 
every Riemannian geodesics on the base $B$ starting at $b$ 
admits a unique lift to the subriemannian
geodesic on $Q$ starting at $q$ with the velocity vector in $\mathcal T^{hor}$.
\end{thm}

\begin{figure}[ht!]
\input epsf
\centerline{\epsfysize=0.33\vsize\epsffile{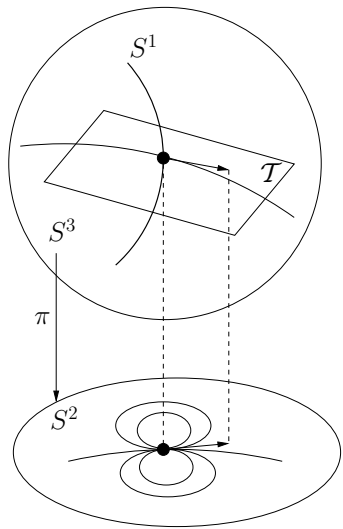}}
 \caption{Projections of subriemannian geodesics from  $(S^3, \mathcal T)$
in the Hopf bundle give circles in $S^2$, 
only one of which, the equator, is a geodesic on the base $S^2$.}\label{fig:hopf}
\end{figure}

\begin{exmp}\label{hopf}
{\rm
Consider the standard Hopf bundle $\pi:S^3\to S^2$, with the two-dimensional
distribution  $\mathcal T$ transversal to the fibers $S^1$. Fix the standard 
metric on the base $S^2$ and lift it to a subriemannian metric on 
$S^3$, which defines a subriemannian submersion. 
If the distribution $\mathcal T$ is orthogonal to the fibers, the manifold
$(S^3, \mathcal T)$ can locally be thought of as 
the Heisenberg 3-dimensional group. 
Then all subriemannian geodesics on $S^3$ with a given horizontal 
velocity  project to a 1-parameter family of 
circles on $S^2$ with a common tangent element. However, 
only one  of these circles, the equator, is a geodesic on the 
standard sphere $S^2$. Thus the equator can be uniquely lifted to a subriemannian
geodesic on $S^3$ with the given initial vector. 

Note that the uniqueness of this lifting holds even if the distribution  
$\mathcal T$ is not orthogonal, but only transversal, say at
a fixed angle, to the fibers $S^1$, see Figure \ref{fig:hopf}. 
}
\end{exmp}

%\begin{figure}[ht!]
%\begin{center}
%\centerline{\scalebox{0.5}{\input{mostranfig5.pstex_t}}}
% \caption{Projections of subriemannian geodesics from  $(S^3, \mathcal T)$
%in the Hopf bundle give circles in $S^2$, 
%only one of which, the equator, is a geodesic on the base $S^2$.}\label{fig:hopf}
%\end{center}
%\end{figure}

\proofbegin of Theorem \ref{thm:sub-reduction}.
To prove this theorem we describe the Hamiltonian setting of the 
subriemannian submersion.

Let $Ver$ be the vertical subbundle in $TQ$ (i.e. tangent planes to the fibers
of the projection $Q\to B$). 
Define $Ver^\perp\subset T^*Q$ to be the corresponding annihilator, i.e.
$Ver^\perp_q$ is the set of all covectors $\alpha_q\in T^*_qQ$ 
at the point $q\in Q$ which annihilate the vertical space $Ver_q$.

\begin{defn}
{\rm
The restriction of the
subriemannian exponential map $\exp^\tau:T^*Q\to Q$ to the
distribution $Ver^\perp$ is called the {\it horizontal
exponential}
$$
\overline{\exp^\tau}: Ver^\perp\to Q
$$
and the corresponding geodesics are the {\it horizontal subriemannian geodesics}.
}
\end{defn}

%Note that the subriemannian exponential map admits such a restriction: 
%a geodesic started with the velocity in $Ver^\perp$ stays there, 
%cf. Proposition \ref{redy}.
The symplectic reduction identifies the quotient
$Ver^\perp/G$ with the cotangent bundle $T^*B$ of the base. 
Note that the subdistribution
$\mathcal T^{hor}$ defines a horizontal bundle
for the principal bundle $Q\to B$ in the usual sense.
The definition of subriemannian submersion (translated to the cotangent spaces,
where we replace  $\mathcal T^{hor}$ by $Ver^\perp$)
gives that the subriemannian Hamiltonian $H^{\mathcal T}$ defined by
(\ref{subHam}) descends to a Riemannian Hamiltonian $H^{B,\mathcal
T}$ on $T^*B$. Moreover, Hamiltonian trajectories of $H^{B,\mathcal
T}$ starting at the cotangent space $T^*_b B$ are in one-to-one
correspondence with the trajectories of $H^{\mathcal T}$ starting
at the space $Ver^\perp_q$. The projection of these Hamiltonian
trajectories to the manifolds $B$ and $Q$ via the cotangent bundle
projections $\pi^{T^*B}$ and $\pi^{T^*Q}$,  respectively, gives
the result. 
\proofend

\begin{cor}
For a subriemannian submersion, geodesics on the base give rise only 
to normal geodesics in the total space.
\end{cor}

In order to describe the geodesic geometry on the tangent, rather than cotangent,
bundle of the manifold $Q$,  we fix a
Riemannian metric on $Q$ whose restriction to the distribution $\tau$ 
is the given subriemannian metric $\langle\,,\,\rangle^\tau$.
This Riemannian metric allows one to identify 
the cotangent bundle $T^*Q$ with the tangent bundle $TQ$. 
Then the exponential map $\exp^\tau$ can be viewed as a map $TQ\to Q$. 
It is convenient to think of $\mathcal T^{hor}$ as the
horizontal bundle 
and identify it with the annihilator $Ver^\perp$.
This way  horizontal subriemannian geodesics are
geodesics with initial (co)vector in the horizontal bundle  $\mathcal T^{hor}$.
This identification is particularly convenient for the
infinite-dimensional setting, where we work with the tangent
bundle of the diffeomorphism group.

\bigskip

%%%%%%%%%%%%%%%%%%%%%%%

\subsection{A subriemannian analog of the Otto calculus.}

Fix a Riemannian metric $\langle\, ,\,\rangle ^M$ on the manifold $M$. 
Let $P^\tau:TM\to\tau$ be
the orthogonal projection of vectors on $M$ onto
the distribution $\tau$ with respect to this  metric. Let
$(\nu,\eta_1)$ and $(\nu,\eta_2)$ be two tangent vectors in the
tangent space at the point $\nu$ of the smooth Wasserstein space.
Recall that for a fixed the volume form $\mu$, we define the
subriemannian Laplacian as $\Delta^\tau f:=\di_\mu (P^\tau \nabla f)$.

Define a {\it nonholonomic Wasserstein metric} as the (weak) Riemannian metric on
the (smooth) Wasserstein space $\W$ given by
\begin{equation}\label{smetric}
\langle (\nu,\eta_1),(\nu,\eta_2)\rangle ^{\W,\T}:=\int_M
\langle P^\tau\nabla f_1(x),P^\tau\nabla f_2(x)\rangle ^M\nu\,,
\end{equation}
where functions $f_1$ and $f_2$ are solutions of the subriemannian Poisson equation
\[
-(\Delta^\tau f_i)\nu=\eta_i\,
\]
for the measure $\nu$.

\begin{thm}\label{submain}
The geodesics on the Wasserstein space $\W$ equipped with the nonholonomic
Wasserstein metric (\ref{smetric}) have the form
$(\overline{\exp^\tau}(tP^\tau\nabla f))_*\nu$,
where $\overline{\exp^\tau}: \mathcal T^{hor}\to M$ is the horizontal
exponential map and $\nu$ is any point of $\W$.
\end{thm}

To prove this theorem we first note that the Riemannian metric
$\langle \,,\,\rangle ^\diff$ defined on the diffeomorphism
group restricts to a subriemannian metric $\langle \,,\,\rangle ^{\diff,\T}$
on the right invariant bundle $\T$.

\begin{prop}
The map $\pi:(\diff,\T)\to\W$ is a subriemannian submersion of the
subriemannian metric  $\langle \,,\,\rangle ^{\diff,\T}$
on the diffeomorphism group with distribution
$\T$ to the nonholonomic Wasserstein metric $\langle ,\rangle ^{\W,\T}$.
\end{prop}

\begin{proof}
This statement can be derived from the Hamiltonian reduction,
similarly to the Riemannian case.

Here we prove it by an explicit computation.
Recall that the map $\pi:\diff\to\W$ is defined by $\pi(\phi)=\phi_*\mu$.
Let $X\circ\phi$ be a tangent vector at the point $\phi$ in the diffeomorphism
group $\diff$. Consider  the flow $\phi_t$ of the vector field $X$, and note that
$\pi(\phi_t\circ\phi)=\phi_{t*}\phi_*\mu$.  To compute the derivative $D\pi$
we differentiate this equation with respect to time $t$ at $t=0$:
$$
D\pi(X\circ\phi)=\LD_{-X}(\phi_*\mu)=-(\di_{\phi_*\mu}X)\phi_*\mu\,,
$$
by the definition of Lie derivative. A vector field $X$ from the
horizontal bundle $\T^{hor}$ has the form
$(P^\tau \nabla f)\circ\phi$, and for it the equation becomes
$$
D\pi((P^\tau\nabla f)\circ\phi)=-(\Delta^\tau f)\,\phi_*\mu\,,
$$
where the Laplacian $\Delta^\tau$ is taken with respect to the volume form
${\phi_*\mu}$.

Therefore, for horizontal tangent vectors
$(P^\tau\nabla f_1)\circ\phi$ and $(P^\tau\nabla f_2)\circ\phi$ at the point $\phi$
their subriemannian inner product is
$$
\langle (P^\tau\nabla f_1)\circ\phi,(P^\tau\nabla f_2)\circ\phi\rangle ^\diff
=\int_M \langle P^\tau\nabla f_1\circ\phi,P^\tau\nabla f_2\circ\phi\rangle ^M\mu\,.
$$
After the change of variables  this becomes
$$
\int_M \langle P^\tau\nabla f_1, P^\tau\nabla f_2\rangle ^M \phi_*\mu
=\langle D\pi((P^\tau\nabla f_1)\circ\phi), D\pi((P^\tau\nabla f_2)\circ\phi)\rangle ^{\W,\T}\,,
$$
which completes the proof.
\end{proof}
\medskip

\proofbegin of Theorem \ref{submain}.
To describe geodesics in the nonholonomic Wasserstein space
we define the Hamiltonian $H^\T:T\diff\to\Real$ by
\begin{equation}\label{sHam}
H^\T(X\circ\phi):=\int_M \langle (P^\tau X)\circ\phi,(P^\tau X)\circ\phi\rangle \mu\,.
\end{equation}
The  Hamiltonian flow with Hamiltonian
$H^\T$, has the form $\exp^\tau ((tP^\tau X)\circ\phi)$ 
according to Theorem \ref{HamFl}.
By taking its restriction to the bundle $\mathcal T^{hor}$ 
and projecting to the base we
obtain that the geodesics on the smooth Wasserstein space are
\[
(\overline{\exp^\tau} ((t P^\tau \nabla f)\circ\phi))_*\nu\,,
\]
where $\nu=\phi_*\mu$ and $P^\tau \nabla f$ 
is defined by the Hodge decomposition for the field $X$.
This completes the proof of Theorem \ref{submain}.
\proofend

\smallskip

\begin{rem}
{\rm
For a horizontal
subriemannian geodesic $\varphi_t(x):=\overline{\exp^\tau}(tP^\tau \nabla f(x))$ 
with a
smooth function $f$, the diffeomorphism  $\varphi_t$ satisfies
$\ddt\varphi_t=(P^\tau \nabla f_t)\circ\varphi_t$ and $f_t$ is the solution
of the Hamilton-Jacobi equation
\begin{equation}\label{subHJ}
\dot f_t+H^\tau(\nabla f_t)=0
\end{equation}
with the initial condition $f_0=f$, see Corollary \ref{HJ}.
This equation determines horizontal subriemannian geodesics on
the diffeomorphism group $\diff$. In the Riemannian case, one can see that
the vector fields $V_t=\ddt\varphi_t=\nabla f_t\circ\varphi_t$ satisfy
the Burgers equation by taking the gradient of the both sides in (\ref{subHJ}),
cf. Proposition \ref{prop:burg}.
Hence Equation (\ref{subHJ}) can be viewed as a subriemannian analog of the
potential Burgers equation in $\diff$.
However, a subriemannian analog of the Burgers equation for nonhorizontal
(i.e. nonpotential) normal geodesics on the diffeomorphism group is not 
so explicit.
}
\end{rem}

\begin{rem}
{\rm
If the function $f$ is smooth, the time-one-map
$\varphi(x):=\overline{\exp^\tau}(P^\tau \nabla f(x))$ along the geodesics
described in Theorem \ref{submain} satisfies the following {\it
nonholonomic} analog of the {\it Monge-Amp\`{e}re equation}:
$h(\varphi(x))\det(D\varphi(x))={g(x)}$,
where $g$ and $h$ are functions on the manifold $M$ defining two densities
$\theta=g\vol$ and $\nu=h\vol$.

Furthermore, for the case of the Heisenberg group this formal solution $\varphi(x)$
coincides  with the optimal map obtained in \cite{AmRi}. The (minus) potential
$-f$ of the corresponding optimal map satisfies
the $c$-concavity condition for $c=d_\tau^2/2$,
where $d_\tau^2$ is the subriemannian distance, cf. Remark \ref{MA1}.
%All the potentials in the optimal maps previously mentioned including
%the ones in Remark \ref{MA1} and Remark \ref{MA2} are all $c$-concave.
}
\end{rem}

\bigskip

%%%%%%%%%%%%%%%%%%%%%%%%%%%%%%%%%%

\subsection{The nonholonomic heat equation}

Consider the heat equation $\partial_t u=\Delta u$ on a function
$u$ on the manifold $M$, where the
operator $\Delta$ is given by $\Delta f=\di_\mu\nabla f$.
Upon multiplying the both sides of the heat equation by the fixed volume form $\mu$,
one can regard it as an evolution equation on the smooth Wasserstein space $\W$.
Note that the right-hand-side of the heat equation gives
a tangent vector $(\Delta u)\mu$ at the point $u\mu$
of the Wasserstein space.
The Boltzmann {\it relative entropy} functional $\Ent:\W\to\Real$ is
defined by the integral
\begin{equation}\label{Boltz}
\Ent(\nu):=\int_M\log(\nu/\mu)\,\nu\,.
\end{equation}
The gradient flow of $\Ent$ on the
Wasserstein space with respect to the metric $\tilde d$ gives the heat equation, see
\cite{Ot}.

Recall that one can define the
subriemannian Laplacian: $\Delta^\tau f:=\di_\mu (P^\tau\nabla f)$
for  a fixed volume form $\mu$ on $M$.
The natural generalization of the heat equation to the nonholonomic
setting is as follows.

\begin{defn}
{\rm
The {\it nonholonomic} (or, {\it subriemannian}) {\it heat equation}
is the equation $\partial_t u=\Delta^\tau u$ on a time-dependent function $u$ on $M$.
}
\end{defn}

Below we show that this equation in the nonholonomic setting also
admits a gradient interpretation on the Wasserstein space.

\begin{thm}
The nonholonomic heat equation $\partial_t u=\Delta^\tau u$ describes
the gradient flow on the Wasserstein space with respect to the
relative entropy functional (\ref{Boltz}) and the nonholonomic Wasserstein
metric (\ref{smetric}).

Namely, for the volume form
$\nu_t:=g_{t*}\mu$ and the gradient $\nabla^{\W,\T}$ with respect to
the metric $\langle \,,\,\rangle ^{\W,\T}$ on the Wasserstein space one has
\[
\pddt\nu_t=-\nabla^{\W,\T} \operatorname{Ent}(\nu_t)=\Delta^\tau(\nu_t/\mu)\mu.
\]
\end{thm}

\begin{proof}
Denote by $(\nu,\eta)$  a tangent vector to the Wasserstein space
$\W$ at a point $\nu\in \W$, where $\eta$ is a volume form of total
integral zero. Let $\Delta^\tau_\nu$ be the subriemannian 
Laplacian with respect to the volume form $\nu$.

Let $h$ and $h_\Ent$ be real-valued functions on the
manifold $M$ such that $-(\Delta^\tau_\nu h)\nu=\eta$ and
$-(\Delta^\tau_\nu h_\Ent)\nu=\nabla^{\W,\T} \Ent(\nu)$ for the entropy
functional $\Ent$.
 Then, by definition of the metric $\langle \,,\,\rangle ^{\W,\T}$ given by (\ref{smetric}),
we have
\begin{equation}\label{non-hol-heat-1}
\langle (\nu,\nabla^{\W,\T}\Ent(\nu)),(\nu,\eta)\rangle ^{\W,\T}
=\int_M \langle P^\tau\nabla h_\Ent(x),P^\tau\nabla h(x)\rangle ^M\nu.
\end{equation}

On the other hand, by definitions of $\Ent$ and the gradient $\nabla^{\W,\T}$ on the
Wasserstein space, one has:
\[
\langle (\nu,\nabla^{\W,\T}\Ent(\nu)),(\nu,\eta)\rangle ^{\W,\T}
:=\ddtz \Ent(\nu+t\eta)
=\ddtz\int_M\Big[\log\Big(\frac{\nu+t\eta}{\mu}\Big)\Big](\nu+t\eta)\,.
\]
After differentiation and simplification the latter expression becomes
$
\int_M\log(\nu/\mu)\,\eta\,,
$
where we used that $\int_M\eta=0$.
This can be rewritten as
\[
\int_M\log(\nu/\mu)\,\,\eta =-\int_M\log(\nu/\mu)\LD_{P^\tau\nabla
h}\nu =\int_M\left(\LD_{P^\tau\nabla h}\log(\nu/\mu)\right)\,\nu\,,
\]
by using the Leibnitz property of the Lie derivative $\LD$ on the Wasserstein space
and the fact that $-(\Delta^\tau_\nu h)\nu=\eta$.
Note that the  Lie derivative  is the inner product with the gradient, and hence
\[
\int\limits_M\left(\LD_{P^\tau\nabla h}\log(\nu/\mu)\right)\,\nu
=\int\limits_M \langle \nabla\log(\nu/\mu),P^\tau\nabla h\rangle ^M\nu
=\int\limits_M \langle P^\tau\nabla\log(\nu/\mu),P^\tau\nabla h\rangle ^M\nu\,.
%=\langle (\nu,\nabla^{\W,\T}\Ent(\nu)),(\nu,\eta)\rangle ^{\W,\T}
\]
Comparing the latter form with (\ref{non-hol-heat-1}), we get
$P^\tau\nabla h_\Ent=P^\tau\nabla\log(\nu/\mu)$, or, after taking the divergence
of both parts and using the definition of function $h_\Ent$,
\[
\nabla^{\W,\T} \Ent(\nu)=-\Delta^\tau_\nu(\log(\nu/\mu))\,\nu\,.
\]

Finally, let us show that the right-hand-side of the above equation coincides with
$-\Delta^\tau_\mu(\nu/\mu)\,\,\mu$. Indeed, the chain rule gives
\[
\LD_{P^\tau\nabla\log(\nu/\mu)}\nu
=\LD_{(\mu/\nu)P^\tau\nabla(\nu/\mu)}\nu
=(\mu/\nu)\LD_{P^\tau\nabla(\nu/\mu)}\nu+d(\mu/\nu)\wedge
i_{P^\tau\nabla(\nu/\mu)}\nu\,.
\]
The last term is equal to
$
(i_{P^\tau\nabla(\nu/\mu)}d(\mu/\nu))
\nu=\LD_{P^\tau\nabla(\nu/\mu)}(\mu/\nu)\,\nu,
$
which implies that
\[
\LD_{P^\tau\nabla\log(\nu/\mu)}\nu=\LD_{P^\tau\nabla(\nu/\mu)}\mu
\]
by the Leibnitz property of Lie derivative.
Thus
$$
\Delta^\tau_\nu(\log(\nu/\mu))\,\nu
=\di_\nu(P^\tau\nabla(\log(\nu/\mu))\nu
=\LD_{P^\tau\nabla\log(\nu/\mu)}\nu=\LD_{P^\tau\nabla(\nu/\mu)}\mu
=\Delta^\tau_\mu(\nu/\mu)\mu\,.
$$

The above shows that the nonholonomic heat equation is
the gradient flow on the Wasserstein space for the same potential
as the classical heat equation, but with respect to 
the nonholonomic Wasserstein metric.
\end{proof}

\bigskip
%%%%%%%%%%%%%%%%%%%%%%%%%%%%%%%%%%%%%%

%%%%%%%%%%%%%%%%%%%%%%%%%%%%%%%%%%%%%%%%

\begin{ackn} We are much indebted to R.~Beals, Ya.~Eliashberg, 
V.~Ivrii, G.~Misiolek,
D.-M.~Nhieu, R.~Ponge, and M.~Shubin for fruitful discussions.
B.K. is grateful to the IHES in Bures-sur-Yvette for its stimulating 
environment. 
This research was partially supported by an NSERC research grant.
\end{ackn}
%%%%%%%%%%%%%%%%%%%%%%%%%%%%%%%%%%%%%%%

\end{document}